\documentclass[a4paper, preprint, review, 11pt]{article}

\usepackage[linesnumbered,ruled,vlined]{algorithm2e}

\SetCommentSty{mycommfont}

\SetKwInput{KwInput}{Input}                
\SetKwInput{KwOutput}{Output}              

\usepackage{caption}
\usepackage{amsmath,amssymb}

\usepackage{adjustbox}
 \usepackage{bbm}
 \usepackage{todonotes}
 \usepackage{graphicx}
\usepackage{color}
\definecolor{cornell-red}{RGB}{179,27,27}
\usepackage{subcaption}
\usepackage{mathtools}
\usepackage{booktabs}
\usepackage{longtable}
\usepackage{dsfont}
\usepackage{float}
\usepackage{natbib}
    \usepackage{multirow}
 
\usepackage[left=1in,right=1in,top=1in,bottom=1in]{geometry}
\usepackage{rotating}
\usepackage{xcolor}
\usepackage[colorlinks = true,linkcolor = blue,urlcolor  = blue,citecolor = blue,anchorcolor = blue]{hyperref}
\usepackage[flushleft]{threeparttable}

\usepackage{booktabs}
\usepackage{multirow}
\usepackage{rotating,tabularx}

\usepackage{xcolor}

\def\mathcolor#1#{\@mathcolor{#1}}
\def\@mathcolor#1#2#3{%
  \protect\leavevmode
  \begingroup
    \color#1{#2}#3%
  \endgroup
}
\usepackage{pifont}
\usepackage{blindtext,mathtools}
\usepackage{geometry}
\usepackage{graphicx}
\usepackage{natbib}
\usepackage{lscape}

 \bibpunct[, ]{(}{)}{,}{a}{}{,}%
 \def\bibfont{\small}%
\begin{document}




\title{Routing mobile health clinics: An integrated routing and resupply plan based on synchronization}


\author{Faisal Alkaabneh$^{1*}$, Sam Jotham Sutharson$^{2}$}
\maketitle

\author{$^{1}$Industrial Engineering, American University of Sharjah, Sharjah, 26666, United Arab Emirates, Email: falkaabneh@aus.edu. $^{*}$\textit{Corresponding Author}}

\author{$^{2}$Industrial \& Systems Engineering, North Carolina A\&T State University, Greensboro, NC 27401, United States, Email: ssutharson@aggies.ncat.edu.}

\abstract{%
As an important means of providing medical services in developing countries and remote areas, Mobile Health Clinics (MHCs) focus on distributing medical supplies and providing basic health needs to underserved communities. In this paper, we propose a new model for the mobile health clinics with resupply from a truck and heterogeneous demand. In addition to the traditional routing decisions, our model also establishes en-route resupply plan. Adding the en-route resupply plan adds complexities to this problem as more constraints need to be added to accommodate for the synchronization between a fleet of MHCs and a resupply truck. We model this problem as a vehicle routing problem with multiple synchronization constraints and heterogeneous demand (VRPMSC-HD) and formulate a mixed integer linear programming (MILP). We propose a metaheuristic based on adaptive large neighborhood search (ALNS) with new operators to solve large-scale instances of the model. To demonstrate the value of the synchronization approach, we compare the total distance traveled by the mobile health clinics and the resupply truck and the arrival of the last mobile health clinic to the depot against a model where mobile health clinics are allowed to perform multiple trips to resupply from the depot. Our results reveal that despite the increase of 6.79\% in traveled distance under the synchronization approach, the reduction in the latest arrival time is 16.06\% on average. Implying that the utilization of the fleet of mobile health clinics can be significantly improved at the cost of extra traveling time when combining the traveling time of all the mobile health clinics and the resupply truck.          
}%

\noindent



%


\subsection*{Highlights}
\begin{itemize}
    \item We present a novel mixed integer programming model that addresses a complex mobile health clinics routing, assignment, and resupply problem.
    \item We design a metaheuristic to solve the challenging NP-Hard model in reasonable computational time with good quality solutions.
    \item Our experimental results highlight the remarkable efficiency of our proposed model, which successfully improves utilization of mobile health clinics by 16\%, on average, by leveraging en-route resupply as opposed to the classical approach where mobile health clinics resupply at the depot. This demonstrates the potential impact of our approach in optimizing patient care management.
\end{itemize}
\textbf{Conflict of Interest:} Authors declare no conflicts of interest.

\textbf{Funding:} This research is partially supported by USDOT Tier 1 University Transportation Center, Center for Freight Transportation for Efficient \& Resilient Supply Chain (FERSC). Any opinions, conclusions or recommendations are those of the authors and do not necessarily reflect the view of the DOT.

\textbf{Data Availability Statement}
The manuscript has associated data in a data repository. 

\section{Introduction}
\label{sec:intro}
Healthcare systems in developing countries face the challenges of very restrictive budget limitations and a growing population that is located in scattered geographic locations. In such a situation, distributing medical supplies and equipment to those populations is a challenge. As a possible way to provide a cost-effective distribution system, some governments and institutions utilize mobile healthcare clinics \cite{mcgowan2020mobile,santa2023multi}. Mobile healthcare clinics (MHCs) improve access to healthcare services and supplies by providing a package of limited primary health services, with referral to nearby fixed structures for conditions not manageable under this package \cite{weinstein2014telemedicine}. Managing a fleet of MHCs requires several decisions, for example, finding the optimal routing of MHCs to visit the selected community locations (see \cite{yucesoy2022mobile}), optimal locations to deploy MHCs (see \cite{li2023integer}), or location-routing where the decision maker decides on the locations of MHCs services, the assignment of patients to selected locations, and the routing of MHCs (see \cite{savacser2022mobile}). \\

We study a problem related to the delivery of medical services and supplies to demand points\footnote{A demand point might be a community living in a remote area, location of under-served population, or simply a site to increase the accessibility to healthcare services} with substantial change to the resupply process to an MHC: instead of returning back to the depot, a truck dedicated to resupply a fleet of MHCs may resupply an MHC at any demand point. We consider the problem from the perspective of a non-profit organization providing these services. Our aim is to improve the operational efficiency of a fleet of MHCs combined with a truck loaded with medical supplies whose only task is to serve as a mobile depot for MHCs by resupplying them with medical supplies once they run out of medical supplies or they are at running low on supplies. The problem involves deciding the locations to be served by each MHC, routing of MHCs, when and where each MHC gets resupplied, and routing of the resupply truck to resupply MHCs. Here we note that deciding the locations where an MHC meets the truck for resupply is an important aspect of our problem, and this is where synchronization comes into picture. The potential locations of a resupply can be at any demand point. We also note that capacity allocation and workforce assignment decisions are not considered in the context of our problem and can be handled as a succeeding decision based on the routes of the mobile facilities.\\

Utilization of vehicle routing with multiple synchronization constraints models to healthcare has recently emerged as a prosperous area, see \cite{polnik2021multistage}. \cite{soares2023synchronisation} state that in complex supply chains, where multiple vehicles, crews, materials and other resources are involved, synchronization can be a catalyst for more efficiently combining different operations that require different resources by decreasing unproductive times. \cite{mankowska2014home} demonstrate the added value of utilizing vehicle routing with synchronization models to optimize routing and scheduling of caregivers to visit patients who require the presence of more than one caregiver simultaneously (or within a specific time) in achieving low traveling cost for caregivers, low average waiting times for patients, a fair distribution of inevitable tardiness, or a combination thereof. Various applications and sectors benefited greatly from using vehicle routing with synchronization models, e.g., perishable items transportation (see \cite{anaya2021iterated}), workforce management (see \cite{goel2013workforce}), agriculture (see \cite{alkaabneh2024matheuristic}), and city logistics (see \cite{nolz2020two}) to name a few. On the other hand, to the best of our knowledge, there is no single paper in the literature that explores the use of vehicle routing with synchronization models in the sector of humanitarian logistics despite the critical role humanitarian logistics play in addressing pressing needs by under-served populations. Hence, the motivation behind this work is to study the effectiveness and efficiency of utilizing vehicle routing with synchronization constraints models to Mobile Health Clinics systems as an example of humanitarian logistics. Namely, we develop a model where a fleet of MHCs works in tandem with a truck for resupplying while en-route without the need to go back to the depot for a resupply. We demonstrate how such model may lead to significant savings in completing service compared to the multi-trip model where vehicles go back to the depot for resupply while on the other hand leads to longer total traveled distance. \\

Motivated by the above observations, this paper investigates improving the operationl efficiency of a fleet of mobile health clinics that serves a heterogeneous demand. Realistic settings are adopted to humanitarian supply chain. A mixed integer programming model is formulated and presented, which determines the optimal routing of MHCs and a resupply truck as well as optimal resupply decisions in terms of schedule and location. Due to the large number of binary variables and the NP-hardness of the developed model, we develop an effective metaheuristic based on Adaptive Large Neighborhood Search (ALNS) to find good quality solutions in reasonable amount of computational time. We conduct extensive numerical experiments to: (i) demonstrate the effectiveness of the proposed metaheuristic and (ii) illustrate the advantages and disadvantages of our proposed model when compared against the classical vehicle routing problem with multi-trips.\\

The remainder of this study is organized as follows. Section \ref{sec:litRev} presents the relevant literature review, focusing on studies on the topic of mobile health clinics and vehicle routing problem with synchronization. In Section \ref{sec:Operation}, we present a clear and detailed description of the proposed model of mobile health clinics with a resupply truck. In Section \ref{sec:mathModel}, we present the mathematical model and in section \ref{sec:ALNS} we outline the metaheuristic we design to solve the problem. In section \ref{sec:Computational} we present the computational results to demonstrate the performance of the developed metaheuristic and we illustrate the advantages and disadvantages of the proposed model as opposed to the classical approach where MHCs can resupply only by going back to the depot. And finally, conclusions and future research directions are presented in Section \ref{sec:conclude}.

\section{Literature Review}
\label{sec:litRev}
A considerable number of papers have been published on the topics of the vehicle routing problem with synchronization or the management of a fleet of mobile health clinics. This paper aims to optimize the delivery services provided by a fleet of MHCs using a vehicle routing model that considers multiple synchronization constraints. To provide a structured overview, we categorize the related literature into three main aspects: (i) the vehicle routing problem with synchronization, (ii) the routing of mobile healthcare facilities, and (iii) the vehicle routing problem with multiple trips.\\

\subsection{Vehicle Routing Problems with Synchronization}
The Vehicle Routing Problem (VRP) with synchronization can be defined as routing problem where routes of vehicles are interdependent \cite{soares2023synchronisation}. Generally speaking, there are more than one type of vehicles needed to perform specific tasks, e.g., more than one nurse need to be present at a patient location who needs more than one specialized operation. Reviews of Vehicle Routing Problem with synchronization are \cite{drexl2012synchronization,soares2023synchronisation}. According to \cite{soares2023synchronisation}, there are four types of synchronization, namely (i) routing with synchronisation of schedules, (ii) routing with transfers or cross-docking requirements, (iii) routing with autonomous vehicles, and (iv) routing with trailers or passive vehicles. Given the features of the problem we study in this paper, our work falls within the third classification. Specifically, we may consider the truck responsible for the resupplying as an active truck while the fleet of MHCs might be viewed as a fleet of passive vehicles. Furthermore, works within the third category can be further classified based on the number of active and passive vehicles, specifically, the model of a single truck and multiple drones working in tandem. In this subsection, we detail, to the best of our knowledge, articles published on the topic of routing a truck and a multi-drone in tandem.\\

The problem of pairing autonomous vehicles (drones) with traditional delivery trucks was introduced by \cite{murray2015flying} as the Flying Sidekick Travelling Salesman Problem (FSTSP). The FSTSP is the problem of finding routes for a truck paired with a single drone to deliver packages to a given set customers. The drone can be launched from any customer location to deliver a package to a different customer, after delivery completion, the drone flies back to meet the truck at a different customer location from which the launch happened. The truck, in the meantime, may visit more than one customer to arrive at the retrieval site. Note that the model of  \cite{murray2015flying} assumes that the drone visits one customer before getting back to the truck. Numerous studies have since explored variations of the single-truck single-drone problem, for example multiple trucks, multiple trucks and multiple drones, multiple drones with adjustable speeds and battery consumption rate as a function of the payload, and a single truck with unlimited number of drones.

\cite{chang2018optimal} study a single-truck multi-drone problems in which the truck deploys multiple drones from distributed launch sites along the truck’s route and the drones return to the truck before the truck departs for its next destination. A done may visit one customer before returning to the truck and the proposed solution methodology was k-clustering with routing.

\cite{karak2019hybrid} consider vehicle routing problem with pickups and deliveries in which multiple drones are mounted on a single vehicle and drones are allowed to visit up to two customers before returning back to the truck. In their model, customers are served only by drones, a truck is only used as a mobile depot for the drones, and drones may be launched and retrieved at a subset of the nodes called `station nodes'.

\cite{poikonen2020multi} introduce and model the k-multi-visit drone routing problem where a truck is equipped with $k$ drones for delivery with heterogeneous packages. In their model, a drone may visit up to two customers before the need to go back to the truck for swaping/recharging batteries or reloading packages. They model the power consumption of drones as a function of drone load. The proposed approach in their work to solve the developed mathematical models was based on flexible heuristic solution that consists of three phases: route, transform, and shortest path. Under the same setups, namely a truck with multiple drones, \cite{leon2022multi} propose an agent-based approach in which locations to be visited are the agents instead of orders or vehicles.

\cite{murray2020multiple} present an extension to the ﬂying sidekick TSP in which a single truck and multiple heterogeneous drones work in tandom for deliver parcels and they modeled the queueing of drone launch and retrieval activities. The solution methodology in their work was based on three-phase approach. In their model, a drone can only visit one customer before re-joining the truck.

\cite{gu2022hierarchical} consider a fleet of trucks, each truck is equipped with a single drone that can perform up to two visits before the need to go back to the truck. Their work models drone endurance as a function of both the payload and flight time rather than distance or flying time only. The solution methodology employed in their work is based on hybrid iterated local search and variable neighborhood search algorithm.

In the field of nonprofit OR, \cite{bonku2024collaborative} study the problem of managing a fleet of mobile food pantries to distribute food supplies to agencies serving local communities. The problem was modeled as a vehicle routing problem with multiple synchronization constraints. The developed models focuses on efficiency (allocate as much food as possible to minimize food waste), equity (ensure fair allocation across agencies), and effectiveness (minimize routing cost)  measures. The authors solved the developed models using a commercial solver using randomly generated instances and present extensive analysis and managerial insights on the structure of the optimal solutions.

It is worth mentioning that are studies in the literature that assume that deliveries can be made only by the drones while the truck's job is to launch and retrieve drones \cite{mathew2015planning,bin2017routing,peng2019hybrid}. However, in all these studies, the drone may visit one customer only before getting back to the truck to reload.

In light of these recent developments, our study contributes to the literature by specifically addressing truck and a fleet of MHC where an MHC may visit any number of nodes as long as the capacity of an MHC is not exceeded with heterogeneous demand of customers. Furthermore, an MHC may get a resupply at any node rather than restricting the resupply operation to a predefined set of nodes.

\subsection{Routing Mobile Healthcare Clinics}
Mobile Healthcare Clinics provide critical services to remote communities or undeserved populations in both developed and underdeveloped countries. In the United State, for example, the St. David's foundation dental program operates a fleet of 9 MHCs to provide access to free dental care for low-income communities in Central Texas \cite{jackson2007creating}. Furthermore, MHCs provide a wide range of services and medical procedures, such as otological \citep{lim2021clinical}, pediatric \citep{dawkins2013dental}, dental \citep{sepulveda2022scheduling}, and optometry \citep{ruggeiro1995evaluation}. In almost all these studies, routing the fleet of MHC, scheduling, and allocation of communities to served nodes are the central part of the decision making process.

One of line of research that is also relevant to routing MHCs is the mobile facility routing. The mobile facility routing problem is concerned with establishing routes for a fleet of mobile facilities to maximize demand covered by these facilities within a planning horizon \cite{halper2011mobile}.

\cite{doerner2007multicriteria} study location-routing problem for a fleet of mobile healthcare facilities. They developed a multiobjective combinatorial optimization with the following objectives: (i) economic efficiency based on to the tour length, (ii) an approximation criteria based distances to the nearest tour stops corresponding to p-median location problem formulation, and (iii) a coverage criterion measuring the percentage of the population unable to reach a service location with a given radius.

\cite{csahinyazan2015selective} design a mobile blood collection system with the primary objective of increasing blood collection levels, collected blood is a perishable product and has to be returned to the depot at the end of each day. A special vehicle called `shuttle' brings back the collected blood to the depot. The authors propose a two-stage integer programming model and a heuristic algorithm to determine the tours of the bloodmobiles and the shuttle. This work is close to our work as there are two types of vehicles; however, three major differences can be found: (1) their model decides on which nodes to visit; our model on the other hand imposes a constraint that all nodes has to be visited, (2) their model decides on daily schedule of bloodmobile and imposes the requirement that donated blood has to be collected by the shuttle in the same day or the bloodmobile has to return to the depot in the same day; in our model, we relax this requirement and allow an MHC to wait to get resupply, and (3) their model assumes that a bloodmobile can visit at most one demand point per day while our formulation uses continuous time space to model the problem allowing for more flexibility.

\cite{yucel2020data} propose a data-driven optimization approach for scheduling and routing of mobile medical services. In particular, they view the problem from the perspective of service provider whose goal is to maximize the revenues generated by bringing the service closer to potential customers. They develop a mixed integer programming model based on a variant of the team-orienteering problem with prizes coming from covered locations.

\cite{sepulveda2022scheduling} study a mobile dental clinic scheduling problem motivated by real-world challenge faced by a school-based mobile dental care program in Southern Chile. They model the problem as a parallel machine scheduling problem with sequence-dependent setup costs and batch due dates with the aim of maximizing equity. They solve the developed model using genetic algorithm.

\cite{santa2023multi} present a multi-period location routing problem for mobile clinics with the aims to select communities to be served and create routes for the mobile clinics. The objective function of their model aims to maximize coverage and continuity of care functions. Their work was motivated by a collaboration with the Premi$\grave{e}$re Urgence Internationale.

The model proposed in this study is flexible and reproducible for any mobile clinic deployments, independent of the service offered, the planning horizon, or the number of products.

\subsection{Vehicle Routing Problem with Multiple Trips}
A variation of the vehicle routing problem where each vehicle performs multiple trips within the planning horizon is known as the Vehicle routing problem with multiple trips (MTVRP). This problem is more representative of a wide range of real-life scenarios \cite{Brandão1998The} and city logistics \cite{Cattaruzza2017Vehicle}. The vehicle routing problems with multiple trips are crucial in modeling and solving real life logistics problems, where vehicles are expected to make more than one trip within the planning horizon \cite{Cattaruzza2016The}.

\cite{Cattaruzza2018Vehicle} provides a comprehensive review of the literature that exists in the domain of vehicle routing problems with multiple trips. The authors define MTVRP with a set of trips and an assignment of a trip to each vehicle where an objective like the travelling time is minimized and the following conditions are met  i) each trip starts and ends at the depot , ii) each customer is visited exactly once, iii) the sum of the demands of  all the customers any one trip does not exceed the capacity of the vehicle, and iv) the sum of the duration of the trips taken by one single vehicle across multiple trips does not exceed the maximum duration or to a single vehicle throughout the planning horizon. It is also pointed out that some of the vehicles may not be used based on the configuration of the problem. The travel times that are associated with multi-trip problems are usually the symmetric matrix, and the problem can be indifferently defined on either a directed or an non directed graph.Several mathematical formulations including exact methods, heuristic approaches, and hybrid meta-heuristic algorithms have been developed to solve different types of vehicle routing problems with multiple trips as discussed in \cite{Pan2020Multi-trip}.

Algorithms like adaptive large neighborhood search (ALNS)  and variable neighborhood descent (VND) have been shown to produce good results in handling the complexities of the multi-trip vehicle routing problem. The objective functions of these problems often aim to minimize the total travel distance or the total time when constrained by parameters such as the vehicle capacities time windows and maximum trip duration 
per vehicle \cite{Cattaruzza2014A}.

\cite{Sahin2022A} incorporates a fleet of heterogeneous vehicles and multiple depots, which is a variant of the MTVRP. The heterogeneous vehicle fleet includes different kinds of vehicles that have different travel times based on the type of the vehicle, and the network consists of multiple depots instead of one. The heterogeneity of vehicles and the presence of multiple depots introduces a layer of complexity to accurately model real world city logistics operations, highlighting the utility of the MTVRP over a traditional VRP.
 
\cite{Battarra2009An} presented an iterative solution approach to reduce the number of multi trips, in turn reducing the number of vehicles that are required by the fleet. The proposed iterative solution relied on the decomposition of the problem into smaller problems each of which was solved by a specific heuristic to produce feasible solutions. The authors refer to this class of problems as the minimum multiple trip VRP ( MMTVRP). The heuristics that were used to solve the decomposed problems were guided by an adaptive guidance approach. The efficiency of the combination of the adaptive guidance approach and decomposition mechanism was computationally verified by applying to a real-world data set. This method successfully reduced the minimum number of vehicles required. 

\cite{pisinger2007general} presented a unified ALNS based heuristic which was able to solve 5 different variations of the vehicle routing problem, and pointed out that this approach had several advantages. One of the key advantages of the method was its ability to adjust the different weights that were assigned to the modifiers, effectively enabling it to self-calibrate to a certain extent. The results showed that this heuristic framework improved on 183 out of 483 benchmark solutions for the vehicle routing problems. 

\cite{Martínez2013A} utilized a two-step heuristic to solve large and real-world size instances of a multi-trip vehicle routing problem with time windows. The two-step approach included a modified Solomon sequential insertion heuristic followed by the Tabu Search. The Tabu search was it likes to improve the initial solution that was the output of the modified sequential insertion heuristic.

Although MTVRP models are designed to capture the complexities of real-world problems compared to those of classical VRP models, several shortcomings do exist and have been identified in literature \cite{Cattaruzza2018Vehicle}. Existing literature does not directly compare a synchronized routing strategy against the MTVRP. The synchronization constraints can theoretically reduce the number of trips made by introducing dependencies between multiple vehicles. This paper directly addresses the shortcomings of the multi-trip models by proposing the routing strategy where the mobile health clinics are synchronized with the single resupply vehicle, eliminating the need for multiple trips. This effectively increases the range of the MHC as the central depot may be located much farther away from the remote nodes. The resupply vehicle effectively acts as a mobile depot to resupply the MHC’s en-route.

\section{Operation description}
\label{sec:Operation}

Mobile health clinics (MHCs) facilitate the delivery of care and medical supplies to areas that need to receive medical assistance for groups of populations that may otherwise have difficulties in obtaining the required medical care and supplies owing to a lack of traditional infrastructure to provide them. In this section, we provide clear and detailed information of the problem we study in this work. The problem is modeled as a vehicle routing problem with multiple synchronization constraints. Using the notation of the vehicle routing literature, a node represent a geographical location that need to be served via a vehicle (i.e., mobile health clinic). Based on the implementation in practice, the nodes can be small towns, villages, or even an emergency medical outpost that serves the locality. An MHC contains equipment and supplies that would either provide primary medical care in situations like a natural disaster or specialized care when there are situations like providing vaccine shots or screening for diseases.\\

The model used in this work assumes that all the MHCs in the fleet are homogeneous and are outfitted to carry out the same function in a given instance. The mobile clinics can carry multiple types of products like vaccines, medicines, and test kits, but all the products together share the same storage space and cannot exceed the total capacity of an MHC. The mobile clinic’s operation cycle starts at the central depot where they are loaded with the required supplies and they visit the nodes to provide care and medical supplies. Once a service is completed at a given node, they move on to the next node. Once the MHC supplies are depleted, the MHC is resupplied using a resupply truck. The resupply truck is primarily a robust transportation truck that lacks the equipment or personel required to serve a node and its only task is to resupply the MHCs with the supplies they need, namely, the resupply truck can not serve a node and it's sole function is to resupply an MHC. Since they only store and transport the supplies, the resupply truck may travel faster than an MHC.\\

\begin{figure}[h!]
  \centering
  \includegraphics[scale = 0.50]{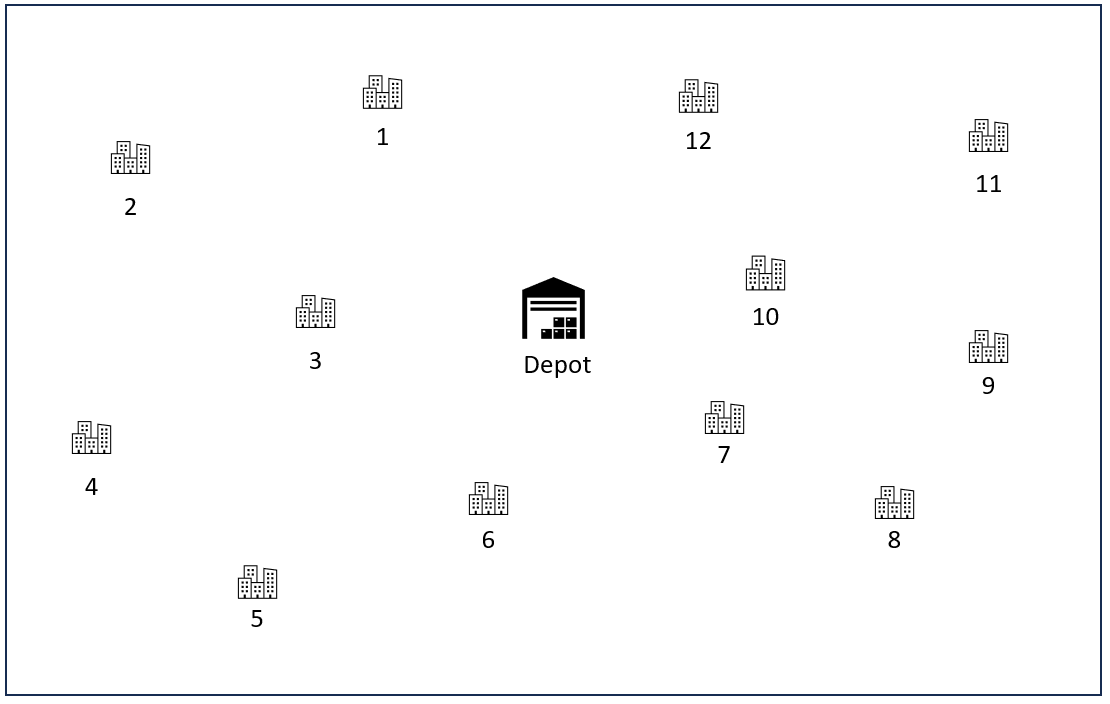}
  \caption{Node and depot locations }
  \label{fig:empty}
\end{figure}

In this model, the quantity of demand that is needed at each node is known
in advance. Similarly, the service time needed at each node also known in advance. Each node requires a specific type of product (or a set of products). The MHC will not travel to a new node if it does not have enough supplies to serve that node. If the demand of the next node cannot be satisfied with the supplies on board, the resupply truck is needed to the current node to provide the supplies to the MHC. If the resupply truck is being used at a different node to resupply another MHC in the fleet or in its way to the current node, the former MHC will have to wait till the resupply truck completes that task and travels to the current node. This is illustrated using a simple example using figures \ref{fig:empty} and ref{fig:filled}. Note that an MHC waiting to get supplies is a waste of time resources and decreases the utilization of an MHC; hence, it is desired to reduce the waiting time of an MHC to resupply.\\

\begin{figure}[h!]
  \centering
  \includegraphics[scale = 0.50]{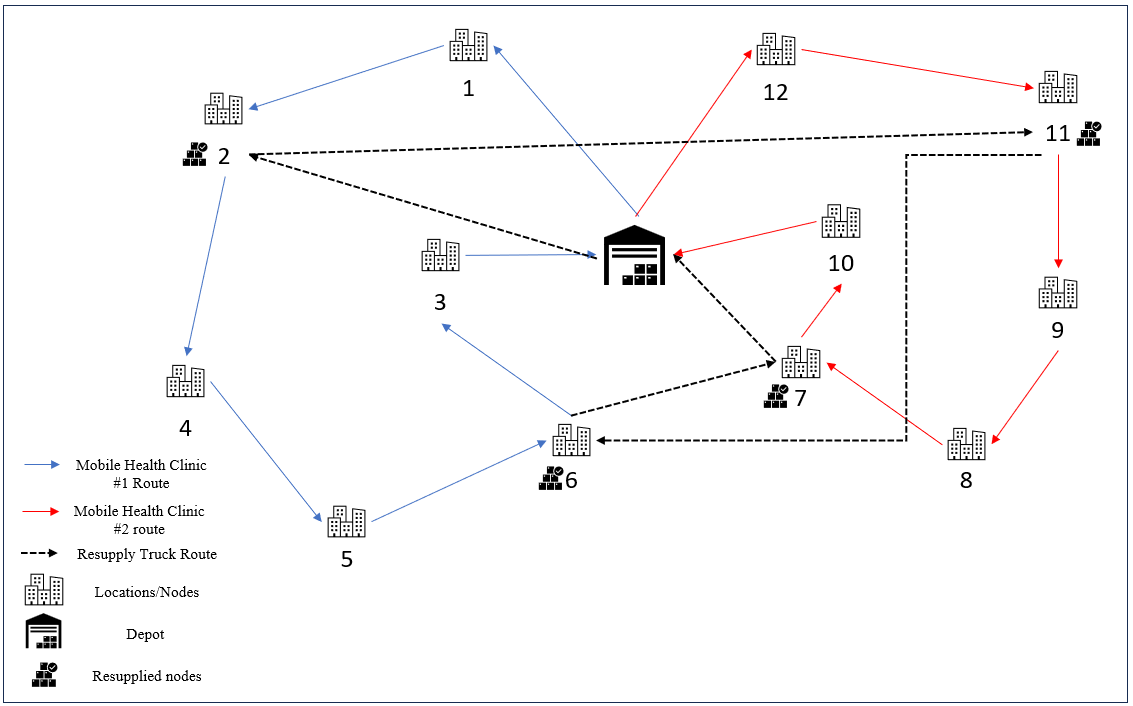}
  \caption{Mobile Health Clinic problem with Resupply Truck }
  \label{fig:filled}
\end{figure}

Figure \ref{fig:empty} shows the 12 locations that are the nodes used in this example. Each location is a node that is numbered from 1 to 12. There is a central depot that is approximately centrally located with respect to all the nodes. The routing for this example model is shown in Figure \ref{fig:filled}. There are two MHCs in this example, and they serve their own routes consisting of nodes 1, 2, 4, 5, 6 and 3 for MHC 1 represented by the blue arrows, and nodes 12, 11, 9, 8, 7 and 10 for MHC 2 represented by the red arrows. The resupply truck routing is represented by the black dotted arrows.\\

In this example, MHC 1 travels to node 1, serves it, and moves on to node 2. After serving node 2, it does not have enough supplies to serve the next node in the route, node 4, and the resupply the truck is summoned. Ideally, the model that we used in this paper routes the resupply vehicles and mobile clinics to minimize the wait time of the MHC as waiting time represent waste of utilization of an MHC. Therefore, in this case, the routing from the model will ensure that the resupply truck is already at location 2 and there will be little to no waiting time as it is the first node in the resupply truck's route. In the meantime, MHC 2 completed serving nodes 12 and 11 and needs to resupply at node 11. In this case, it is likely that MHC will wait the resupply truck at node 11 till the resupply truck arrives to start the resupply process. Once resupply process is completed, MHC 2 can move to serve nodes 9, 8, an 7 before it runs out of supplies again. Note that MHC 1 started serving nodes 4, 5, and 6 simultaneously and before leaving node 6, it needs a refill. Afer the resupply truck completed resupplying MHC 2 at node 11, it should have started moving to node 6 in anticipation of the resupply processes needed for MHC 1 after serving node 6. Note that it is permissible for the resupply truck to wait for the MHC to finish service before the resupply process starts. Once MHC 1 resupplied at node 6, it moves to node 4 and then finally back to the depot. On the other hand, the resupply truck will move to node 7 to resupply MHC 2. Finally, the resupply truck goes back to the depot while MHC 2 continues serving node 10 and then returns back to the depot. \\

It is important to emphasize that MHCs are critical and expensive resources as means to provide medical support. Hence, it is important to minimize the total routing, waiting, and service times of a fleet of MHC in a \textit{balanced way}. Balanced way refers to the fact that workload of one MHC within a fleet of MHC should not significantly exceed the workload of other MHCs within the same fleet \cite{matl2018workload}.

\section{Mobile Health Clinic with Resupply Truck Synchronized Routing Problem and Heterogeneous Demand}
\label{sec:mathModel}
\subsection{Mathematical model}
\label{sec:model}

The vehicle routing problem with multiple synchronization constraints and heterogeneous demand can be characterized as:
\begin{itemize}
    \item a set of nodes with known demand,
    \item a set of MHCs with limited capacity responsible for supplying demand nodes with products and services,
    \item a single resupplying truck responsible for resupplying the set of MHCs, 
    \item the resupplying of an MHC may only happen at a demand node and only after the completion of a service at a node. This is a mild assumption we add to simplify the model, in practice a resupply may happen before or after a service completion, 
    \item it is permissible that an MHC or the resupplying truck wait at the refilling location, and 
    \item each demand node shall receive its full demand. 
\end{itemize}

The problem considers deciding which MHC visits which node, the nodes of resupplying the MHCs, the routing of MHCs, and the routing of the supply truck. First, in section \ref{sec:notation} we introduce the notation and then in Section \ref{mathmodel} we provide the complete mathematical model. 

\subsubsection{Mathematical notation}
\label{sec:notation}
The set of demand nodes is denoted as $\mathcal{N}_c$, the depot is considered node $0$ and the set of all nodes in the network is denoted as $\mathcal{N}=\lbrace0\rbrace\cup \mathcal{N}_c$. We assume that the network is fully connected and hence there is an arc from node $i\in\mathcal{N}$ to node $j\in\mathcal{N}:i\neq j$, the traveling time from node $i\in\mathcal{N}$ to node $j\in\mathcal{N}$ is denoted as $t_{ij}$ and it is proportianal to the distance between the nodes. On the other hand, the traveling time from node $i\in\mathcal{N}$ to node $j\in\mathcal{N}$ for the supply truck is denoted as $r_{ij}$. The MHC fleet size is denoted as $numMC$ and we assume that MHCs are homogeneous in terms of capacity and traveling speed, in case of having non-homogeneous fleet of MHCs, an index for the set of MHC shall be introduced. The total capacity of an MHC is denoted as $Q^{MC}$. The set of products is denoted as $\mathcal{K}$ and the demand of each product at each node is known and denoted as $d_{ik}$. The replenishment time needed to resupply an MHC is denoted by $\xi$ and we assume that the replenishment time is constant and does not depend on the number of units available at the MHC at the time of replenishment. A service time at each demand node $i\in\mathcal{N}_c$ is known and denoted as $s_i$. And finally we define $M$ as a sufficiently large number.\\

As of decision variable, let $x_{ij}$ be a binary decision variable that takes the value of 1 if an MHC travels from node $i\in\mathcal{N}$ to node $j\in\mathcal{N}:i\neq j$ and 0 otherwise. Let $y_{ij}$ be a binary decision variable that takes the value of 1 if the resupply truck travels from node $i\in\mathcal{N}$ to node $j\in\mathcal{N}:i\neq j$ and 0 otherwise. Let $z_i$ be binary decision variable that takes the value 1 if resupply happens at node $i\in\mathcal{N}_c$ and 0 otherwise. Let $l$ be a positive continuous decision variable denoting the latest arrival of an MHC to the depot. Let $c_i$ be a positive continuous decision variable denoting the service completion time of an MHC at node $i\in\mathcal{N}_c$. Let $w_i$ be a positive continuous decision variable denoting the waiting time of an MHC at node $i\in\mathcal{N}_c$ and let $u_i$ be a positive continuous decision variable denoting the waiting time of the resupply vehicle at node $i\in\mathcal{N}_c$. Let $a_i$ be a positive continuous decision variable denoting the arrival time of the resupply vehicle at node $i\in\mathcal{N}_c$. And lastly, let $q_{ik}$ be a positive continuous decision variable denoting the quantity of product $k\in\mathcal{K}$ received by an MHC at node $i\in\mathcal{N}_c$ from the resupply truck and let $g_{ik}$ be a positive continuous decision variable inventory level of product $k\in\mathcal{K}$ available by an MHC at at the arrival at node $i\in\mathcal{N}_c$.

Table~\ref{Table:Notation} summarizes all notation. Using this notation, we formulate the MILP presented in the next subsection:
\begin{table}[h!]
\small
  \renewcommand{\arraystretch}{0.4}
\caption{Notation}
\label{Table:Notation}
\begin{tabular}{ll}
\hline 
\multicolumn{2}{l}{\textbf{Sets}}                       \\ \hline
$\mathcal{N}$         & Set of nodes                        \\
$\mathcal{N}_c$         & Set of nodes to be served                        \\
$\mathcal{K}$         & Set of products                        \\
\hline
\multicolumn{2}{l}{\textbf{Indices}}                       \\ \hline
$i,j$         & Node indices                        \\
$k$           & Product                        \\ \hline
\multicolumn{2}{l}{\textbf{Parameters}}                        \\ \hline
$numSp$           & Number of sprayers  \\
$d_{ik}$          & Number of units of product $k\in\mathcal{K}$ to be supplied at node $i\in\mathcal{N}_c$          \\
$s_i$             & Service time of node $i\in\mathcal{N}_c$          \\
$t_{ij}$    & Travelling time between nodes $i$ and $j$ for an MHC      \\
$r_{ij}$    & Travelling time between nodes $i$ and $j$ the resupply truck     \\
$Q_{MC}$   & Capacity of a MHC \\
$\xi$      & Time needed to resupply a mobile health clinic \\
$\gamma$      & Time needed to refill the tanker at the depot \\
$M$       & Big positive number \\ \hline
\multicolumn{2}{l}{\textbf{Decision variables}} \\ \hline
$x_{ij} $ & Binary variable that is equal to 1 if and only if MHC travels from node $i\in\mathcal{N}$ to node $j\in\mathcal{N}$ \\
$y_{ij}$     & Binary variable that is equal to 1 if and only if the resupply truck travels from node $i\in\mathcal{N}$ to node $j\in\mathcal{N}$ \\
$z_i$     & Binary variable that is equal to 1 if and only if there is a resupply at node $i\in\mathcal{N}_c$ \\
$a_{i}$     & Arrival time of the resupply truck at node $i\in\mathcal{N}_c$\\
$c_i$     & Time of arrival of an MHC at node $i\in\mathcal{N}_c$ \\
$g_{ik}$     & Number of units of product $k\in\mathcal{K}$ available at the MHC upon arrival to node $i\in\mathcal{N}_c$ \\
$q_{ik}$     & Number of units of product $k\in\mathcal{K}$ resupplied to an MHC at node $i\in\mathcal{N}_c$ \\
$l$       & arrival time of the latest MHC to the depot  \\
$u_i$       & Waiting time of an MHC at node $i\in\mathcal{N}_c$ \\
$w_i$       & Waiting time of the resupply truck at node $i\in\mathcal{N}_c$ \\
 \hline
\end{tabular}
\end{table}

\subsubsection{Mathematical model formulation}
\label{sec:mathmodel}

\begin{eqnarray}
\label{eq:FObjt}
\displaystyle{\min}&& \displaystyle{l+\sum_{i\in\mathcal{N}}\sum_{j\in\mathcal{N}:i\neq j}r_{ij}y_{ij}}\\
s.t.
\label{eq:MCRouting1}
& &\displaystyle{\sum_{i\in \mathcal{N}:i\neq j} x_{ij} = 1 \quad \forall \quad j\in \mathcal{N}_c,} \\
\label{eq:MCRouting2}
& & \displaystyle{\sum_{i\in\mathcal{N}_c:i\neq j} x_{ij} = \sum_{i\in\mathcal{N}_c:i\neq j} x_{ji} \quad \forall \quad } j\in \mathcal{N}, \\
\label{eq:MCRouting3}
& & \displaystyle{\sum_{i\in\mathcal{N}_c} x_{0i} =\sum_{i\in\mathcal{N}_c} x_{i0} = numMC,} \\
\label{eq:RSRouting1}
& & \displaystyle{\sum_{i\in\mathcal{N}_c:i\neq j} y_{ij} = \sum_{i\in\mathcal{N}_f:i\neq j} y_{ji} \quad \forall \quad  j\in \mathcal{N},}\\
\label{eq:RSRouting2}
& & \displaystyle{z_{i} = \sum_{j\in \mathcal{N}:j\neq i}y_{ji} \quad \forall \quad i\in\mathcal{N}_c,}  \\
\label{eq:RSRouting3}
& & \displaystyle{\sum_{j\in \mathcal{N}_c}y_{0j}\leq 1,}  \\
\label{eq:MCcomp1}
&&\displaystyle{c_{j} \leq c_{i}+w_i+s_j+t_{ij}x_{ij}+\xi z_i + M(1-x_{ij}) \quad \forall i \in \mathcal{N},j \in \mathcal{N}_c: i\neq j,}\\
\label{eq:MCcomp2}
&&\displaystyle{c_{j} \geq c_{i}+w_i+s_j+t_{ij} x_{ij}+\xi z_i - M(1-x_{ij}) \quad \forall i \in \mathcal{N},j \in \mathcal{N}_c: i\neq j,}\\
\label{eq:ResuComp1}
& & \displaystyle{a_{j} \geq a_{i}+u_{i}+\xi z_i+r_{ij}y_{ij} - M(1-y_{ij}) \quad \forall i \in \mathcal{N},j \in \mathcal{N}_c: i\neq j,}\\
\label{eq:ResuComp2}
& & \displaystyle{a_{j} \leq a_{i}+u_{i}+z_i\omega+r_{ij}y_{ij} + M(1-y_{ij}) \qquad \forall i \in \mathcal{N},j \in \mathcal{N}_c: i\neq j,}\\
\label{eq:waiting1}
& & \displaystyle{w_{i} \geq a_{i} - c_{i} \quad \forall i \in \mathcal{N}_c,}\\
\label{eq:waiting2}
& & \displaystyle{u_{i} \geq c_{i} - a_{i} \quad \forall i \in \mathcal{N}_c,}  \\
\label{eq:latestArrival}
& & \displaystyle{l \geq c_{i} + t_{i0} x_{i0}\quad \forall i \in \mathcal{N}_c,}\\
\label{eq:inv1}
& & \displaystyle{g_{jk} \leq g_{ik}-d_{ik}x_{ij}+q_{ik} + Q^{MC}(1-x_{ij}) \quad \forall i \in \mathcal{N},j \in \mathcal{N}, k \in \mathcal{K}: i\neq j,}\\
\label{eq:inv2}
& & \displaystyle{g_{jk} \geq g_{ik}-d_{ik}x_{ij}+q_{ik} - Q^{MC}(1-x_{ij}) \quad \forall i \in \mathcal{N},j \in \mathcal{N}, k \in \mathcal{K}: i\neq j,}  \\
\label{eq:inv3}
& & \displaystyle{g_{jk} \geq d_{jk}\sum_{i\in\mathcal{N}:i\neq j}x_{ij} \quad \forall j \in \mathcal{N}_c,j \in \mathcal{N}, k \in \mathcal{K},}\\
\label{eq:inv4}
& & \displaystyle{\sum_{k \in \mathcal{K}}g_{ik} \leq Q^{MC} \quad \forall i \in \mathcal{N}_c,}\\
\label{eq:refill1}
& & \displaystyle{\sum_{k \in \mathcal{K}}q_{ik} \leq Q^{MC}z_i\quad \forall i \in \mathcal{N}_c,}\\
\label{eq:FIntg}
& &\displaystyle{y_{ij}, x_{ij},z_{i} \in \lbrace 0,1\rbrace}, s_i,u_i,w_i,a_{i},q_{ik},l,g_{ik} \geq 0 \quad \forall i,j,k. \
\end{eqnarray}

The objective function (\ref{eq:FObjt}) aims at minimizing the latest arrival of an MHC back to the depot after the completion of its route plus the total routing time of the resupply truck. The reasoning behind the selection of this objective function is to provide an equitable working load for the MHC drivers as discussed in \cite{matl2018workload}. The decision to minimize the latest arrival over distance also ensures that all nodes are served within a maximum prescribed time, prioritizing patient care over distance constraints. \cite{tan2023multi} showed that time based objective function ensures that the model adapts well to uncertainties of real world conditions. \cite{liu2023systematic} also highlights the importance of minimizing completion times, leading to proper utilization of resources as a key advantage. On the other hand, for the resupply truck, the focus is on minimizing the total traveled distance.\\

Constraints (\ref{eq:MCRouting1}) ensure that each demand node is visited once. Constraints (\ref{eq:MCRouting2}) are the inflow and outflow constraints implying that if an MHC visits a node, it shall leave that node. Constraints (\ref{eq:MCRouting3}) imply that the number of MHC used in service equals exactly the number of available MHCs. Constraints (\ref{eq:RSRouting1}) guarantee that if the resupply truck enters a node, it leaves that node. Constraints (\ref{eq:RSRouting2}) state that if a resupply happens at node $i\in\mathcal{N}_c$, the resupply truck has to visit that node. Constraints (\ref{eq:RSRouting3}) imply that only a resupply truck is available by restricting the number of arcs leaving the depot to be 1 at most. \\

Constraints (\ref{eq:MCcomp1})-(\ref{eq:MCcomp2}) calculate the arrival time of an MHC at node $j\in\mathcal{N}_c$ after service completion at node $i\in\mathcal{N}_c$. Note that the arrival time equals to the traveling time plus the service completion time at the prior node if no resupply was needed or traveling time plus the end of the resupply time at the prior node if resupply was needed. Likewise, constraints (\ref{eq:ResuComp1}) and (\ref{eq:ResuComp2}) calculate the arrival time of the resupply truck at the successive node based on the waiting time to start the resupply process, the starting time for the resupply, and the traveling time from node $i\in\mathcal{N}$ to node $j\in\mathcal{N}$. Constraints (\ref{eq:waiting1}) and (\ref{eq:waiting2}) calculate the waiting time of an MHC and the resupply truck, respectively. Constraints (\ref{eq:latestArrival}) calculate the latest arrival of an MHC back to the depot. Constraints (\ref{eq:inv1}) and (\ref{eq:inv2}) keep track of the number of units available at an MHC after service completion at a demand node and in case of a resupply. Constraints (\ref{eq:inv3}) guarantee that the available number of units needed for any product with an MHC is at least the demand of that product at a demand node. Constraints (\ref{eq:inv4}) respect the holding capacity of an MHC. Constraints (\ref{eq:refill1}) are the linking constraint between the quantity resupplied for an MHC and the resupply process itself. Lastly, constraints (\ref{eq:FIntg}) are the domain constraints of the decision variables.\\

Mathematical model (\ref{eq:FObjt})-(\ref{eq:FIntg}) is at least an NP-hard model as it is a generalization of the classical VRP model with synchronization constraints. To solve large-scale instances of the developed model, we advise a metaheuristics based on Adaptive Large Neighborhood Search (ALNS).

\section{Metaheuristic for Vehicle Routing Problem with Multiple Synchronization Constraints and Heterogeneous Demand}
\label{sec:ALNS}
The problem we discuss in this study is at least an NP-hard problem since it is a VRP with multiple synchronization constraints. The mathematical formulation presented in Section \ref{sec:mathmodel} can be solved using broad-sense solvers such as Gurobi only when the size of the instances is small. Given the need to solve large-scale instances in the real-world settings, the problem may be intractable. As a result, we develop a metaheuristic to find high-quality solutions for large instances within reasonable practical computational times. The metaheuristic we develop consists of two phases, the first phase aims at finding an initial feasible solution that respects all constraints and the second phase aims at improving the initial feasible solution iteratively using Adaptive Large Neighborhood Search (ALNS).\\

ALNS have been successful in solving larger instances of combinatorial problems. The neighborhood search is a good candidate to solve complex combinatorial routing problems as demonstrated in previous research (e.g., \cite{mourad2021integrating,alkaabneh2023routing}).\\

The fundamentals of the ALNS method are based on the principle of destroying and recreating a solution to explore more feasible solutions in the neighborhood of the existing solution. The implementation of ALNS in this paper involves four steps, namely partial destruction of a given solution, reconstruction of a solution, new solution evaluation, and updating of weights of the destroy and repair operator. The framework used in this paper starts off by generating an initial seed solution, this seed solution is taken as the input for the ALNS method. The generation of the seed solution follows a greedy algorithm, and this is discussed further in Section \ref{sec:feasible solution}. The generated seed solution is then passed to the destroy operation. The destroy operation performs the function of partially removing parts of the solution based on a specific criteria. Both the destroyed solution and the parts that were removed from them are then sent on to the repair operation. The repair operation performs the function of rebuilding a new solution from the destroyed solution and the parts of the solution that were removed, and the newly constructed solution is now passed on to the evaluation phase. The evaluation method evaluates the quality of the solution and compares it with the existing global best solution found so far during the search. Section \ref{sec:soldRep} talks in detail about the development and implementation of this method. \\

In the absence of synchronization between the resupply truck and an MHC, the metaheuristic will only aim at finding routes for the fleet of MHCs. However, finding the optimal routing of the truck for a given set of resupplying nodes is also an optimization problem by itself. In our work, the developed metaheuristic uses an algorithm to find the resupply nodes for a given routes of MHCs. This algorithm takes the generated MHCs routing solution as an input and outputs the location of resupply nodes based on the demand information, the capacity of an MHC, and the mobile clinic routing information. This algorithm is an important part of the metaheuristic as the resupply nodes are calculated for every solution that is being estimated, and in theory, this also means that the metaheuristic is also optimizing the resupply nodes for a given routing plan of the mobile health clinic. \\

The remaining sections of the metaheuristic are described as follows: Section \ref{sec:feasible solution} details the procedure of generating an initial feasible solution, Section \ref{sec:soldRep} outlines the process of establishing the truck routing and scheduling of the resupply process for given routing plan of the MHCs, Sections \ref{sec:Destroy} and \ref{sec:Repair} present the destroy and repair operators utilized in our ALNS respectively. Section \ref{sec:acceptCriteria} details the acceptance criteria of our ALNS metaheuristic. Section \ref{sec:weightUp} describes the updating procedure of the weights and finally, we present the pseudo code in section \ref{sec:pseudo code}.

\subsection{Initial feasible solution procedure}
\label{sec:feasible solution}
\begin{figure}[h!]
  \centering
  \begin{subfigure}[b]{0.3\textwidth}
    \includegraphics[width=\textwidth]{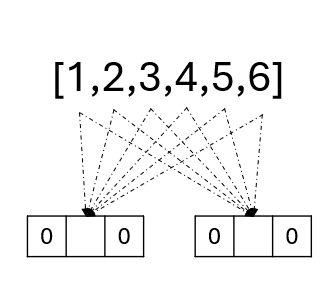}
    \caption{Evaluation of node positions}
    \label{fig:sol-start}
  \end{subfigure}
  \hfill 
  \begin{subfigure}[b]{0.3\textwidth}
    \includegraphics[width=\textwidth]{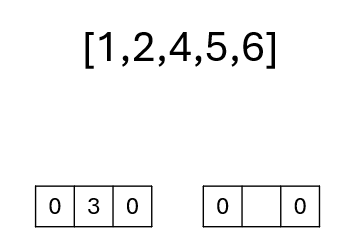}
    \caption{Best node inserted into route and deleted from available node list}
    \label{fig:greedy-insert}
  \end{subfigure}
  \hfill 
  \begin{subfigure}[b]{0.3\textwidth}
    \includegraphics[width=\textwidth]{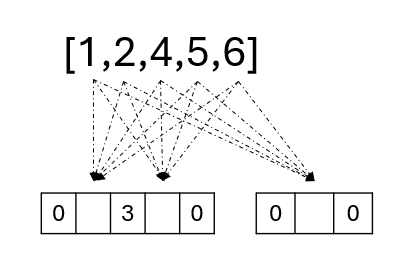}
    \caption{Evaluation of subsequent nodes in all possible positions}
    \label{fig:route-building}
  \end{subfigure}
  
  \caption{Initial solution procedure}
  \label{fig: initial solution}
\end{figure}

Generating an initial feasible solution starts by producing empty routes that equal the number of mobile health clinics that are present in the instance, and each empty route consists of the starting and the ending nodes of these mobile health clinics (namely the depot represented by ‘$0$’). The next step involves estimating the position cost of adding a node to be served by an MHC. For each node that has not been inserted yet, the algorithm estimates the cost of inserting each node at each candidate location and calculates the cost of MHC routes using objective function (\ref{eq:FObjt}) while ignoring the synchronization aspect. At each iteration, the algorithm inserts the node that yields the lowest objective function value at the best position, removes that nodes from the list of uninserted nodes, and repeats these steps till all nodes have been inserted, as shown in Figure \ref{fig: initial solution}. \\

Figure \ref{fig: initial solution} displays how the seed solution is generated for an instance with six nodes and two MHCs. In the represented example, all the available nodes are inserted sequentially into the two available positions in the sub-routes.

\subsection{Calculating resupply locations, waiting times, and routing of the resupply truck}
\label{sec:soldRep}

We design a simple heuristic to calculate the synchronization between MHCs and the truck for given MHC routing plen. The developed heuristic leads to local optimal solutions but it can be excecuted very quickly. \\

We use two sets of sequences to represent a solution of the VRPMSCHD, the first set of sequences represent the routing of MHCs, denoted as $\mathcal{R}= \{ \mathcal{R}_1,\mathcal{R}_2,...,\mathcal{R}_{nM}\}$ where each $\mathcal{R}_m$ is the ordered set of nodes visited by MHC $m$ and the second sequence is the ordered sequence of nodes visited by the resupply truck denoted as $\mathcal{P}$.\\

The first part of heuristic calculates the points at which resupply an MHC shall happen to prevent an MHC running out of products by tracking the quantity remaining of each product $k\in\mathcal{K}$ for each MHC after leaving each node.  In route $\mathcal{R}$, resupply point are calculated by tracking the node at which the vehicle runs out of products to deliver to the next node in the route. The remaining inventory of product $k$ at node $i_n$ is given by
\begin{equation}
g_{i_n,k} = g_{i_{n-1},k}-d_{i_{n-1},k}
\end{equation}
If the remaining inventory of product $k$ at an MHC at node $i_n$ is less than the demand $d_{i_{n+1}k}$ of node $i_{n+1}$, then node $i_n$ is added to the resupply list of nodes. The set of resupply nodes is now represented by the set $\mathcal{P}$ . This is now passed to the second part of the heuristic.\\

The second part of the heuristic calculates the actual arrival of MHC to the nodes and the actual time of resupply starts (i.e., the synchronization). The heuristic starts by picking one of the MHCs and calculates the arrival time to each node $i$ within the selected route until a node was preceded by a resupply node (i.e., node $i-1$ is a resupply node), in this case the heuristic records the arrival time of node $i-1$ in a tuple list with the following format $\mathcal{T}=\lbrace (i, c_i), (j, c_j),....,\rbrace$. This process repeats till all routes have been checked. The next step is to decide the sequence of visiting nodes in list $\mathcal{T}$, the order of refilling then follows the ascending order of nodes based on the values of $c_i$. \\

The third part of the heuristic calculates the waiting times, arrival of the truck, arrival of the MHCs to the next nodes that are preceded by resupply nodes. This process repeats (the second and the third parts of the heuristic) till no resupply node is not synchronized for the arrival of the truck and an MHC.

A roulette-wheel system controls the selection of the destroy and repair operators of a neighborhood as established earlier. All the destroy and repair operators have an equal chance/probability of being chosen at the initial iteration. After a fixed number of iterations, the probabilities of the destroy operators are updated based on their performance.\\ 

\subsection{Destroy operators}
\label{sec:Destroy}
A total of $9$ destroy operators were used in the ALNS method.

\subsubsection{Random Node Destroy}

This destroy operator randomly selects a predetermined number of nodes from a single route in the given input route. The number of nodes to be removed is between  10 and 15 percent of the total number of nodes in the route.  \\

\subsubsection{Longest Node Cost Destroy Operator}

This destroy operator calculates for each node the total distance traveled by an MHC by adding the distance between the previous node and the next node for a given route. Then, this list of node distances is sorted and the top $10\%$ to $15\%$ percent of the nodes in this list are then removed from their routes and inserted in the removed nodes.\\
        
\subsubsection{Destroy Resupply Nodes}

This destroys operator encourages the exploration of the solution space by completely removing the nodes associated with resupply. For a given list of resupply nodes, this operator randomly selects half of the nodes in the list, removes them from their routes, and adds these nodes to the list of removed nodes. 
    
\subsubsection{Destroy Nodes with Wait Times}

This destroy operator removes all the nodes that have a wait time associated with them. 
    
\subsubsection{Destroy Entire Route}

This destroy operator selects an entire route from the main route and removes all the nodes from the selected route and adds the selected nodes to the removed nodes list.

\subsubsection{Destroy Longest Route}

This destroy operator selects the MHC responsible for the latest arrival, removes all the nodes of this route, and adds the selected nodes to the removed list of nodes. 

\subsubsection{Destroy Nodes After Resupply Nodes}

This destroy operator removes the nodes that follow a resupply node. The operator selects 50 percent of resupply nodes at random, and they remove all the nodes that follow the chosen refill nodes till the next resupply node or the depot is encountered. By targeting these specific sets of nodes, the algorithm will be able to exploit the solution space of the previous solution. 

\subsubsection{Destroy Nodes Prior Resupply Nodes}

This destroy operator removes the nodes that precede a resupply node. The operator selects 50 percent of resupply nodes at random, and they remove all the nodes that precede the chosen refill nodes till the next resupply node or the depot is encountered. By targeting these specific sets of nodes, the algorithm will be able to exploit the solution space of the previous solution. 

\subsubsection{Historical Destroy operator}

The historical destroy operator maintains a record of the best objective function for each node at each position. For the sample instance, this operator maintains a table containing the objective function value of each node for each position over every iteration of the ALNS method. When this operator is chosen, the nodes and positions that have the worst values for the objective function are chosen and the top $10\%$ to $15\%$ of them are destroyed from the route.

\subsection{Repair Methods}
\label{sec:Repair}
Three repair operators were used to repair the destroyed routes for the ALNS algorithm. The repair operators were not assigned weights but were chosen at random. The repair operators take the input from the destroy operators output. As discussed in \ref{sec:Destroy}, the destroy operators output a destroyed solution and the nodes that were removed from the original solution. The repair operators used different types of greedy reconstruction algorithms to produce or 'repair' the destroyed solution to give the new solution that will be used by the estimation algorithm to estimate their objective function value. 

\subsubsection{Greedy Repair}

The greedy repair operator uses a greedy algorithm to repair the routes with the destroyed nodes. This operator takes the nodes one by one and places them in all possible positions in the destroyed route to calculate the objective function value. This is then tabulated as a table and then sorted so that the node with the best objective function value and the corresponding position is at the top. This node is then inserted into the route and eliminated from the list of removed nodes. This process is then repeated over and over again till the whole route is completed and the number of routes in the 

\subsubsection{Regret Repair}
The regret repair operator inserts the destroyed nodes into the solution by measuring the regret value. The regret value is the difference between two positions of insertion of the destroyed nodes that produce the best solution. A regret repair of degree 1 calculates the regret value between the top two positions and creates a list that is sorted by the values. The node is inserted into the position with the highest regret value and the remaining destroyed nodes are inserted into the solution in the same way till all the destroyed nodes are inserted into the solution. The third repair operator is a regret repair operator of degree $2$, where the regret value is calculated between the best position and the third best position, skipping the second one, and the process is repeated till all destroyed nodes are replaced.

\subsection{Acceptance Criteria}
\label{sec:acceptCriteria}
The acceptance criteria that are used to accept a new solution as the best solution, to be next current solution, or discard the new solution are a combination of deterministic improvement and simulated annealing. On the completion of an iteration, if the new solution produces an objective function value that is better than the previous global best solution it is accepted as the new global best solution as well as the new solution will be the current solution for the next iteration to be destroyed and repaired. If the new solution does not produce a better objective function value compared to the current solution, it can still be accepted as the new solution based on the probabilistic acceptance using simulated annealing. The probability of acceptance is determined by the difference in solution cost and the temperature of the metaheuristic as shown in \ref{eq:21}.\\
\begin{equation}
    e^\frac{f_{new} - f_{current}}{temperature}
    \label{eq:21}
\end{equation}
where the temperature of the metaheuristic is prescribed to gradually decrease over time, $f_{new}$ is the objective function value of the new solution and $f_{current}$ is the objective function value of the current solution. The value of temperature decreases over iterations ensures that worse solutions have a lesser chance of being accepted with more iterations. This dual acceptance criterion helps the ALNS method in avoiding getting stuck at local minima when searching for the global optimum solution.

\subsection{Updating the Weights}
\label{sec:weightUp}

The weight updating part of the ALNS stores the weights of the nine destroy operators discussed in \ref{sec:Destroy} and updates them once every $150$ iterations. The weights of the nine destroy operations are tracked by $wd_{d}$ where $d$ is the index of the destroy operator used. A performance score $ps$ is used to update the weight of each operator based on its performance in the prior 150 iterations. If the operator contributed to bring a new global solution, the value of the performance score $ps_d$ is set to 10. If the operator contributed to bring a new solution that is better than the current solution, the value of the performance score $ps_d$ is set to 7. If the operator did not contribute to bring a new global solution nor a better new solution but the new generated solution was accepted following the acceptance criteria described in Section \ref{sec:acceptCriteria}, the value of the performance score $ps_d$ is set to 5. And finally, the lowest score $ps_d$ is assigned to the operator that failed to produce a better solution and also failed the acceptance criteria. \\ 

A constant $\gamma\in (0,1)$ is used as the factor that controls the balance between the performance of a destroy operator and the existing weight. The weight updating is carried out by the following equation,
\begin{equation}
    wd'_d=\gamma \cdot wd_d + (1-\gamma) \cdot ps_d
\end{equation}
where $wd'_d$ is the new weight assigned to the $d^{th}$ destroy operator during after every $150$ iterations. Once the new weights are calculated for all the ten operators, they are then normalized and now serve as the probability distribution for selection of the destroy operators.\\

The values of vector $ps_d$ are reset to $0$ after $150$ iterations when the weights of the operators are successfully updated.  This periodic update encourages the algorithm to pick and favor destroy operators that consistently produce good solutions. This does not completely eliminate the methods that have less success in producing solutions and thus maintains diversity in the search process.

\subsection{ALNS Framework: pseudo code}
\label{sec:pseudo code}
Algorithm \ref{ALNSalgo} outlines the complete ALNS metaheuristic developed in this study.

 \begin{algorithm}[h!]
  \renewcommand{\arraystretch}{0.5}
  \small
\DontPrintSemicolon
  \KwInput{A feasible solution $sol_0$ produced by a construction heuristic and a set of ALNS parameters}
  \KwOutput{The best known solution $sol^*$}
  Set $iter \gets 0,\ sol^* \gets s_0,\ sol_{current} \gets sol_0 $, $\mathbf{wd} = [1,...,1]$\;
  \While{$iter < iterMax$}
   {
   		Select a destroy operator $des$ and repair operator $rep$ by a roulette wheel mechanism using weights $\mathbf{wd}$\;
   		Apply destroy operator $des$ on $sol_{current}$ and get $sol^\prime$\;
   		Apply repair operator $rep$ on $sol^\prime$ and get $sol_{new}$\;
   		\If{$f_{sol_{new}} < f_{sol^*}$}
    {
        $sol^* \gets sol_{new},\ sol_{current} \gets sol_{new}$\;
        \label{line8}
        $NoImproveCounter \gets 0$\;
        
    }
    \Else
    {     $NoImproveCounter \gets NoImproveCounter + 1$\;
    	  \If{$f_{sol_{new}}  < f_{sol_{current}}$}
    {
        $sol \gets sol_{new}, sol_{current} \gets sol_{new}$\;
        \label{line13}
        
    }
    \Else
    {
    	$probAccept \gets e^{\frac{f_{sol_{new}}-f_{sol_{current}})}{Temp}}$\;
    	Generate a random number $U\sim Uniform[0,1]$\;
    	\If{$U < probAccept$}
            {
                $sol_{current} \gets sol_{new}$    \tcp*{Accept the new solution which is worse than the current solution for diversification}

            }
            \Else
            {
            	$sol_{current} \gets sol_{current}$ \tcp*{Reject the new solution and use the current solution}
            	
            }
    }
    }
    $iter\gets iter+1$\;
    \If{$iter \% segmentLength == 0$}
        {Calculate $\mathbf{ps}$ and update $\mathbf{wd}$ and\;}

    \If{$NoImproveCounter \% MaxNoImprov == 0$}
        {$sol \gets sol^*,\ NoImproveCounter \gets 0$}
   }
\caption{ALNS framework for the VRPMSC-MT}\label{ALNSalgo}
\end{algorithm}

\section{Computational Results}
\label{sec:Computational}

In this section we conduct computational experiments on test instances using Solomon's instances\footnote{https://www.sintef.no/projectweb/top/vrptw/} to study the performance of the developed metaheuristic and quantify the added value of using the synchronization approach as opposed to multi-trip model where vehicles return back to the depot for a resupply. We begin by introducing the characteristics of the instances and the parameters we use. We then present the results of using the developed ALNS metaheuristic on solving large-scale instances of the developed model on networks with size of 100 nodes. Lastly, we present managerial insights about the difference between the model with synchronization against the model where MHCs have to return back to the depot for a resupply as opposed to having a dedicated truck for en-route resupply.

The model and ALNS met-heuristic presented in this paper are coded in Python and we use Gurobi 11.0.0 optimization software as the MILP solver. All the experiments were conducted on a computer with Intel(R) Core(TM) i7-8700 CPU @ 3.20GHz processor with 16.0GB RAM.\\

\subsection{Test Instances}
We adapt VRPTW Solomon instances, specifically, we use the customers and depot coordinates while we generate other parameters. The computational studies performed on these instances provide insights into the relationships between the solution quality and the computational efficiency of the proposed algorithms, as well as the structural properties of the instances.

\begin{itemize}
\item Locations of the depot and nodes: we use the three types of networks in Solomon's instances. Namely, Random ``R'', Clustered ``C'', and Random and Clustered ``RC''. 
\item Number of nodes to be served $\vert \mathcal{N}\vert$: 15-25 for small instances, 30-50 for medium instances, and 75 and 100 nodes for large-scale instances. 
\item Number of MHC: $\lbrace 2,3,4,5,6,7,8\rbrace$.
\item Number of products $\vert\mathcal{K} \vert=\lbrace 2,3 \rbrace $.
\item Resupply time: 10 units of time.
\item Capacity of an MHC: 26 units.
\item Demand of each node: randomly generated number that is either 4 or 5 product. For the type of product, randomly generated using a uniform distribution.
\item Service time at each node: 20 units of time
\end{itemize}

\subsection{Computational Performance of the Metaheuristic}
In this section, we report on the computational experiment results obtained by our ALNS metaheuristic. The focuses of this section are first to demonstrate the performance of the developed ALNS metaheuristic in obtaining good-quality solutions in a short computational time and secondly to illustrate the added value of using the designed destroy operators.\\

Table \ref{tab:smallInst} displays the results of solving the model using Gurobi solver and the developed ALNS metaheuristic. The first three columns show the number of nodes, the number of MHC, and the number of products, respectively. The fourth column shows the upper bound value found by Gurobil solver after 4 hours (14,000 seconds) of solution time and the fifth column shows the lower bound value at termination of the solver. The sixth column shows the optimality gap of Gurobi solution, the optimality gasp is calculated as $\frac{UB-LB}{UB}*100\%$. We observe that the optimality gaps are very high with an average of 83.29\% and never below 70\%. This observation highlights the difficulty of this problem. We also believe that the poor value of the lower bound is due to the use of \textit{big M} in the formulation. The seventh and eighth columns show the objective function value of the best solution returned by the ALNS and the computational time in seconds, respectively. Clearly, we observe that the best solution found by the ALNS is better than the best solution found by the solver in all cases but the instance of 15 nodes, 3 MHC, and 3 products. On the other hand, the computational time of the ALNS is far better than the computational time of the solver, 4 hours versus less than 20 minutes.\\

\begin{figure}
    \centering
    \includegraphics[width=1\linewidth]{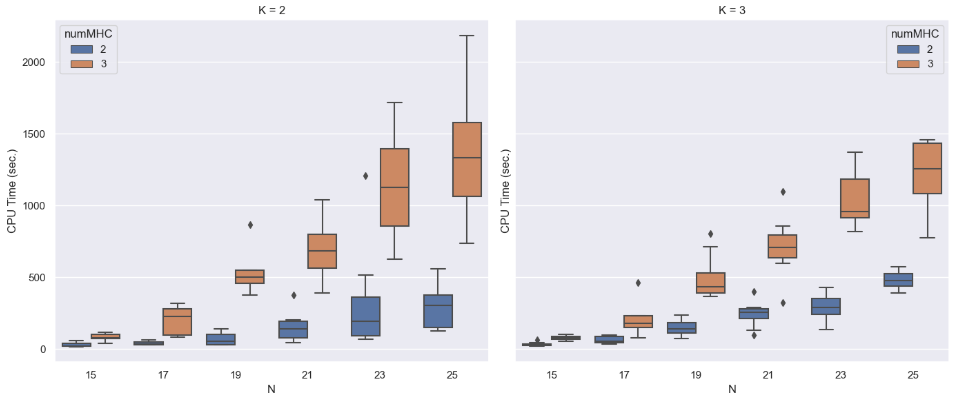}
    \caption{Box plot of computational time for given instance size}
    \label{fig:ALNSBoxPlot}
\end{figure}

To further demonstrate the robustness of the developed ALNS, we generate random instances of 15-25 nodes, under 2 and 3 products, and 2 and 3 MHC. Each instances is solved 10 times under different x-y coordinates of nodes' locations and demand values. We report on the computational time in seconds and the results are summarized in Figure \ref{fig:ALNSBoxPlot}. We observe that the performance of the ALNS is consistent when it comes to the computational time without huge variation between instances of the same number of nodes, products, and number of MHC.

\begin{table}[h!]
\centering
\caption{Computational results of small instances}
\label{tab:smallInst}
\begin{tabular}{cccccccc}
\multicolumn{1}{l}{Nodes} & \multicolumn{1}{l}{numMHC} & \multicolumn{1}{l}{$\vert\mathcal{K}\vert$} & \multicolumn{1}{l}{UB} & \multicolumn{1}{l}{LB} & \multicolumn{1}{l}{Gap (\%)} & \multicolumn{1}{l}{ALNS Obj.} & \multicolumn{1}{l}{ALNS Time (sec.)} \\ \hline
15                        & 2                          & 2                                           & 1007.79                & 221.68                 & 78.13                        & \textbf{993.13}               & 16.39                                \\
15                        & 2                          & 3                                           & 1090.86                & 239.35                 & 78.20                        & \textbf{1081.41}              & 27.18                                \\
15                        & 3                          & 2                                           & 743.61                 & 202.36                 & 72.97                        & \textbf{719.94}               & 41.05                                \\
15                        & 3                          & 3                                           & \textbf{780.58}        & 228.06                 & 70.97                        & 793.72                        & 102.58                               \\
17                        & 2                          & 2                                           & 1107.16                & 204.90                 & 81.61                        & \textbf{1069.90}              & 28.60                                \\
17                        & 2                          & 3                                           & 1222.40                & 226.64                 & 81.58                        & \textbf{1197.48}              & 56.64                                \\
17                        & 3                          & 2                                           & 869.07                 & 161.10                 & 81.59                        & \textbf{827.64}               & 88.63                                \\
17                        & 3                          & 3                                           & 874.27                 & 190.47                 & 78.36                        & \textbf{827.06}               & 232.79                               \\
19                        & 2                          & 2                                           & 1361.58                & 222.22                 & 83.77                        & \textbf{1284.32}              & 32.87                                \\
19                        & 2                          & 3                                           & 1295.87                & 205.62                 & 84.23                        & \textbf{1227.12}              & 74.63                                \\
19                        & 3                          & 2                                           & 1045.94                & 186.88                 & 82.25                        & \textbf{939.49}               & 391.43                               \\
19                        & 3                          & 3                                           & 953.81                 & 166.28                 & 82.69                        & \textbf{922.61}               & 711.63                               \\
21                        & 2                          & 2                                           & 1358.21                & 184.13                 & 86.53                        & \textbf{1274.08}              & 111.48                               \\
21                        & 2                          & 3                                           & 1701.26                & 212.25                 & 87.60                        & \textbf{1500.81}              & 398.55                               \\
21                        & 3                          & 2                                           & 1128.82                & 194.80                 & 82.86                        & \textbf{1072.52}              & 806.18                               \\
21                        & 3                          & 3                                           & 1176.61                & 207.64                 & 82.47                        & \textbf{1096.18}              & 682.68                               \\
23                        & 2                          & 2                                           & 1642.62                & 194.90                 & 88.21                        & \textbf{1536.31}              & 69.50                                \\
23                        & 2                          & 3                                           & 1574.26                & 194.70                 & 87.71                        & \textbf{1464.64}              & 252.94                               \\
23                        & 3                          & 2                                           & 1283.43                & 186.39                 & 85.58                        & \textbf{1182.37}              & 1428.21                              \\
23                        & 3                          & 3                                           & 1256.98                & 163.32                 & 87.10                        & \textbf{1125.59}              & 1372.07                              \\
25                        & 2                          & 2                                           & 1827.61                & 199.83                 & 89.14                        & \textbf{1644.45}              & 123.61                               \\
25                        & 2                          & 3                                           & 1814.06                & 195.90                 & 89.27                        & \textbf{1655.98}              & 467.72                               \\
25                        & 3                          & 2                                           & 1558.96                & 202.90                 & 87.08                        & \textbf{1379.50}              & 1213.45                              \\
25                        & 3                          & 3                                           & 1523.14                & 196.95                 & 87.15                        & \textbf{1283.43}              & 1170.00                              \\ \hline
Avg.                        &                           &                                            &                 &                  & 83.21                        &               & 412.53                               \\
Min.                        &                           &                                            &                 &                  & 89.27                        &               & 1428.21                              \\ 
Min.                        &                           &                                            &                 &                  & 70.97                        &               & 16.39                              \\ \hline
\end{tabular}
\end{table}

\subsubsection{Analysis of destroy operators}

Table \ref{tab:destOp} reports the contribution of each destroy operator in finding a new best solution. More specifically, for all new solutions that resulted in finding a new global best solution, what is the percentage of these solutions resulted from using a given destroy operator. The format of Table \ref{tab:destOp} is as follows, the first column shows the information about an instance in the format of number of nodes-number of MHC-number of products and the rest of the columns are labeled with the given destroy operator. We opt to aggregate all the results across the three types of networks (R, C, and RC) for brevity. \\

\begin{figure}[h!]
    \centering
    \includegraphics[width=0.75\linewidth]{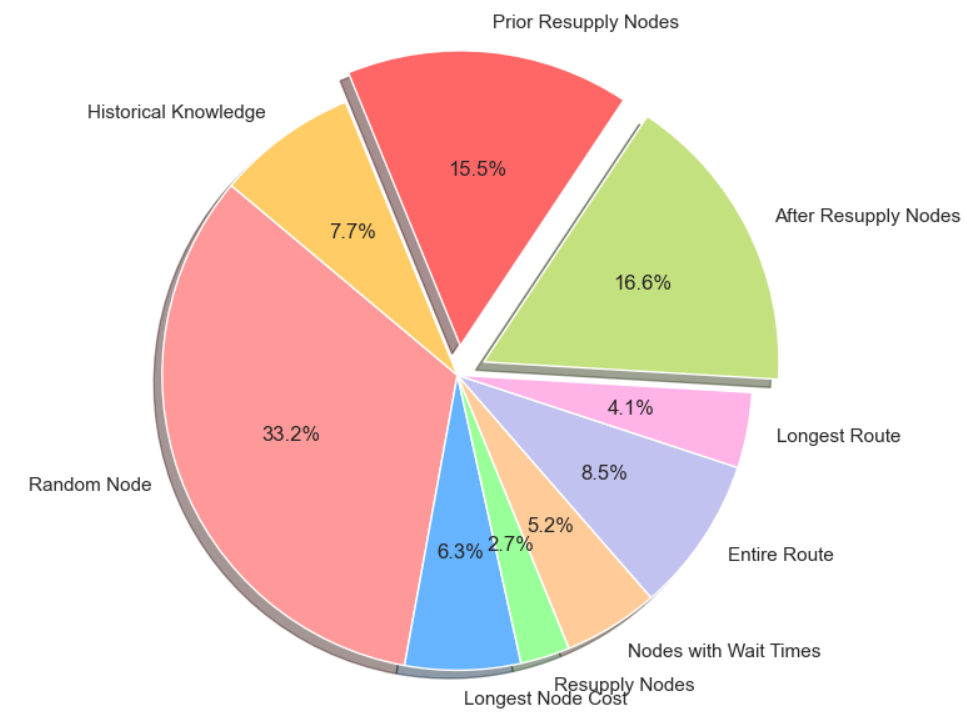}
    \caption{Percentage of getting a new best solution for each destroy operator.}
    \label{fig:destroyPercentage}
\end{figure}

We observe that, on average, the operator that performed the best is the ``Random Node'', followed by the ``Prior Resupply nodes'', and ``After Resupply Nodes'' operators. Interestingly, the operator ``Resupply Nodes'' did not perform well. We believe that the operators of destroy nodes prior to and after supply nodes performed better is due to the large number of nodes typically removed from the solution compared to the operator of ``Supply Nodes'' destroy operator that removes half of the supply nodes. Other operators borrowed from the literature perform well too but not significantly high. Figure \ref{fig:destroyPercentage} display the results as a pie chart. 

\begin{table}[t!]
\caption{Results of effectiveness of each destroy operator}
\label{tab:destOp}
\begin{tabular}{lp{1.5cm}p{1.5cm}p{1.5cm}p{1.5cm}p{1.5cm}p{1.5cm}p{1.5cm}p{1.5cm}p{1.5cm}}
Instance & Random Node & Longest Node Cost & Resupply Nodes & Nodes with Wait Times & Entire Route & Longest Route & After Resupply Nodes & Prior Resupply Nodes & Historical Knowledge \\ \hline
I-30-2-2 & 17.78\%     & 6.67\%            & 2.22\%         & 2.22\%                & 7.78\%       & 6.67\%        & 11.67\%              & 14.44\%              & 11.11\%         \\
I-30-3-2 & 24.71\%     & 4.88\%            & 1.48\%         & 11.45\%               & 9.19\%       & 4.84\%        & 19.25\%              & 17.62\%              & 6.58\%          \\
I-30-4-2 & 25.88\%     & 0.57\%            & 7.18\%         & 5.46\%                & 8.17\%       & 5.07\%        & 22.24\%              & 16.42\%              & 9.00\%          \\
I-30-2-3 & 21.40\%     & 15.20\%           & 4.34\%         & 10.96\%               & 8.57\%       & 7.10\%        & 11.05\%              & 12.97\%              & 8.41\%          \\
I-30-3-3 & 28.01\%     & 4.59\%            & 0.79\%         & 8.76\%                & 11.97\%      & 4.74\%        & 18.04\%              & 17.48\%              & 5.61\%          \\
I-30-4-3 & 30.02\%     & 1.86\%            & 3.51\%         & 8.73\%                & 9.81\%       & 6.37\%        & 15.97\%              & 20.77\%              & 2.96\%          \\
I-40-2-2 & 45.88\%     & 11.37\%           & 0.52\%         & 3.34\%                & 3.85\%       & 7.21\%        & 14.24\%              & 7.94\%               & 5.66\%          \\
I-40-3-2 & 28.65\%     & 5.06\%            & 3.95\%         & 4.21\%                & 7.29\%       & 4.15\%        & 20.18\%              & 17.49\%              & 9.02\%          \\
I-40-4-2 & 34.70\%     & 3.59\%            & 2.54\%         & 4.83\%                & 6.66\%       & 2.20\%        & 25.39\%              & 16.51\%              & 3.58\%          \\
I-40-2-3 & 39.78\%     & 8.87\%            & 2.22\%         & 4.55\%                & 11.38\%      & 1.45\%        & 15.34\%              & 8.45\%               & 7.96\%          \\
I-40-3-3 & 29.53\%     & 3.04\%            & 2.43\%         & 3.63\%                & 10.87\%      & 3.48\%        & 18.24\%              & 19.74\%              & 9.05\%          \\
I-40-4-3 & 26.09\%     & 6.83\%            & 4.40\%         & 3.38\%                & 10.87\%      & 3.84\%        & 20.09\%              & 19.39\%              & 5.11\%          \\
I-50-2-2 & 32.31\%     & 9.81\%            & 0.88\%         & 2.92\%                & 2.33\%       & 1.75\%        & 19.31\%              & 15.48\%              & 15.20\%         \\
I-50-3-2 & 50.99\%     & 3.51\%            & 0.00\%         & 3.51\%                & 7.02\%       & 5.23\%        & 11.36\%              & 11.48\%              & 6.90\%          \\
I-50-4-2 & 38.75\%     & 1.77\%            & 5.22\%         & 4.15\%                & 12.57\%      & 4.09\%        & 15.86\%              & 12.76\%              & 4.82\%          \\
I-50-2-3 & 40.38\%     & 9.72\%            & 4.14\%         & 1.55\%                & 0.76\%       & 2.99\%        & 8.34\%               & 12.77\%              & 19.36\%         \\
I-50-3-3 & 39.68\%     & 9.65\%            & 0.00\%         & 2.73\%                & 6.87\%       & 1.26\%        & 17.14\%              & 19.52\%              & 3.15\%          \\
I-50-4-3 & 37.43\%     & 5.71\%            & 2.45\%         & 6.70\%                & 15.92\%      & 1.22\%        & 11.76\%              & 15.20\%              & 3.60\%          \\ \hline
Avg.     & 32.89\%     & 6.26\%            & 2.68\%         & 5.17\%                & 8.44\%       & 4.09\%        & 16.42\%              & 15.36\%              & 7.62\%          \\
Min.     & 17.78\%     & 0.57\%            & 0.00\%         & 1.55\%                & 0.76\%       & 1.22\%        & 8.34\%               & 7.94\%               & 2.96\%          \\
Max      & 50.99\%     & 15.20\%           & 7.18\%         & 11.45\%               & 15.92\%      & 7.21\%        & 25.39\%              & 20.77\%              & 19.36\%      \\  \hline
\end{tabular}
\end{table}

\subsubsection{Results on large-scale instances}
\begin{figure}[t!]
    \centering
    \includegraphics[width=0.85\linewidth]{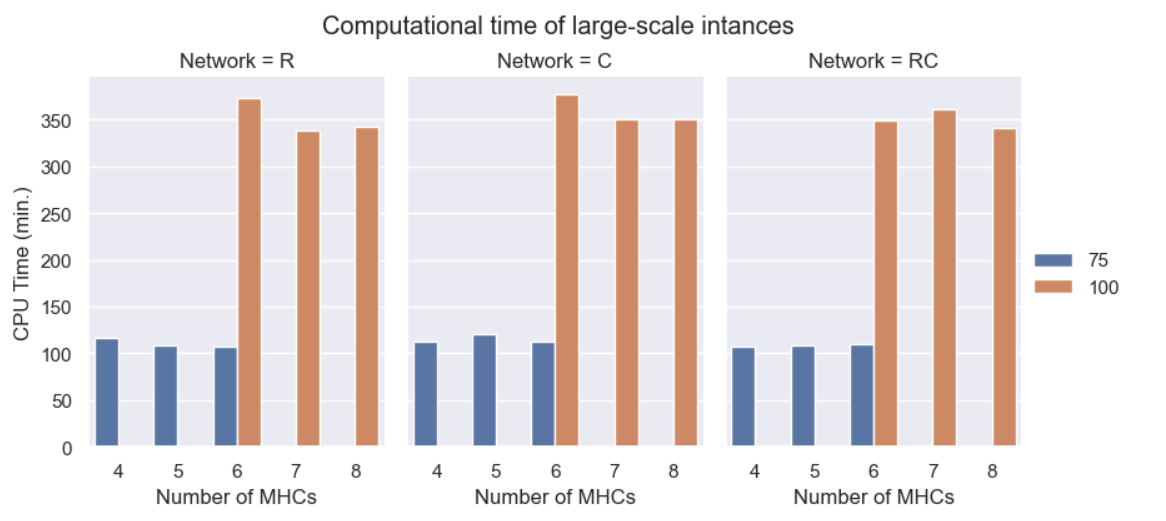}
    \caption{Computational time in minutes for large-scale instances.}
    \label{fig:LSCPUTime}
\end{figure}

In this section, we report on the computational time and the objective function values of running the ALNS metaheuristic on large-scale instances with 75 and 100 nodes; number of MHC of 4, 5, 6, 7, and 8; two products; and networks of R, C, and RC. Results are displayed in Table \ref{tab:largeScale} and Figure \ref{fig:LSCPUTime}.

\begin{table}[t!]
\centering
\caption{Results of ALNS on large-scale instances}
\label{tab:largeScale}
\begin{tabular}{llllp{1.5cm}}
Network & Nodes & $numMHC$ & ALNS\_OBJ & CPU Time (min.) \\ \hline
R       & 75    & 4        & 2920.8    & 116.7           \\
R       & 75    & 5        & 2675.6    & 108.5           \\
R       & 75    & 6        & 2803.9    & 106.8           \\
R       & 100   & 6        & 3560.5    & 373.2           \\
R       & 100   & 7        & 3546.8    & 338.8           \\
R       & 100   & 8        & 3652.6    & 342.2           \\
C       & 75    & 4        & 2804.0    & 112.7           \\
C       & 75    & 5        & 2270.5    & 120.7           \\
C       & 75    & 6        & 2689.6    & 112.1           \\
C       & 100   & 6        & 4196.6    & 376.5           \\
C       & 100   & 7        & 3360.7    & 350.5           \\
C       & 100   & 8        & 3549.8    & 350.8           \\
RC      & 75    & 4        & 3122.7    & 107.4           \\
RC      & 75    & 5        & 3447.7    & 108.3           \\
RC      & 75    & 6        & 3480.8    & 110.0           \\
RC      & 100   & 6        & 4519.9    & 349.2           \\
RC      & 100   & 7        & 4927.4    & 361.1           \\
RC      & 100   & 8        & 4203.4    & 340.5           \\ \hline
Avg.    &       &          &           & 232.6           \\
Min.    &       &          &           & 106.8           \\
Max.    &       &          &           & 376.5    \\ \hline      
\end{tabular}
\end{table}

\subsection{Comparison between models}
The goal of this section is to quantify the advantages and disadvantages of the synchronization model. To this end, we implement a model where MHCs are allowed to go back to the depot for a resupply to ensure that all demand is satisfied. Therefore, we believe that these two models are comparable as they have the same requirements, namely, all nodes are served and resupply is performed. We focus on two performance measures to compare the two models: (1) the total routing distance, and (2) the latest arrival of an MHC. To have a fair and valid comparison, we solve the model developed by \cite{huang2021multi} using Gurobi solver. We selected the model of \cite{huang2021multi} due to the fact that their model takes into account the queuing of the vehicles at the depot for reloading. This is an important factor to consider since the synchronization approach takes into account the waiting time of the resupply truck and the MHCs. The objective function in both models is now set to be minimization of the latest arrival of an MHC.\\

We compare the two models assuming networks of 30 and 35 nodes; three types of networks, R, C, and RC; one product; size of MHC fleet of 3 and 4; capacity of MHC as \textbf{T}ight (22 units) and \textbf{E}xcess (32); and demand of nodes equals to 5 units.\\

Results are displayed in Table \ref{tab:modelsComparison}. The first column indicates the type of network, the second column indicates the capacity of an MHC, the third column indicates the number of MHCs, and the fourth column shows the number of nodes. The fifth column reports the total distance traveled by the fleet of MHCs and the resupply truck. The sixth column reports the arrival time of the latest MHC to the depot under the synchronization model. The seventh column reports the total distance traveled by the fleet of the MHC under the multi-trip model and the eighth column reports the values of the arrival of the latest MHC to the depot. The ninth and tenth  columns show the percentage difference between the two models in terms of traveled distance and arrival of the last MHC, respectively. Note that column 9 is calculated as $\frac{TD_{Synch}-TD_{MT}}{TD_{MT}}*100\%$ where $TD_{MT}$ is the total distance under the multi-trip model and $TD_{Synch}$ is the total traveled distance under the synchronization model. On the other hand, column 10 is calculated as $\frac{LA_{MT}-LA_{Synch}}{LA_{MT}}*100\%$ where $LA_{MT}$ is the arrival time of the last MHC to the depot under the multi-trip model and $LA_{Synch}$ is the arrival time of the last MHC to the depot under the synchronization model.

\begin{table}[t!]
\caption{Results of traveled distance and arrival of the last MHC under the synchronization and multi-trip models}
\label{tab:modelsComparison}
\begin{tabular}{llllllllp{1.5cm}p{1.5cm}}
Netw. & Cap. & $numMHC$ & N  & $TD_{Synch}$ & $LA_{Synch}$ & $TD_{MT}$ & $LA_{MT}$ & Distance Gap (\%) & Arrival Gap (\%) \\ \hline
R       & T        & 3        & 30 & 876          & 465          & 777       & 554       & 12.74             & 16.04            \\
R       & T        & 4        & 30 & 822          & 399          & 806       & 446       & 1.98              & 10.55            \\
R       & E        & 3        & 30 & 807          & 410          & 699       & 536       & 15.39             & 23.50            \\
R       & E        & 4        & 30 & 729          & 338          & 718       & 430       & 1.53              & 21.36            \\
C       & T        & 3        & 30 & 583          & 447          & 530       & 482       & 9.91              & 7.39             \\
C       & T        & 4        & 30 & 563          & 330          & 551       & 387       & 2.25              & 14.71            \\
C       & E        & 3        & 30 & 504          & 348          & 472       & 416       & 6.82              & 16.29            \\
C       & E        & 4        & 30 & 520          & 293          & 507       & 365       & 2.56              & 19.73            \\
RC      & T        & 3        & 30 & 1249         & 575          & 1141      & 677       & 9.50              & 15.02            \\
RC      & T        & 4        & 30 & 990          & 380          & 959       & 411       & 3.23              & 7.66             \\
RC      & E        & 3        & 30 & 833          & 422          & 804       & 497       & 3.61              & 15.09            \\
RC      & E        & 4        & 30 & 893          & 364          & 839       & 511       & 6.44              & 28.72            \\
R       & T        & 3        & 35 & 1030         & 604          & 910       & 688       & 13.22             & 12.16            \\
R       & T        & 4        & 35 & 1157         & 500          & 1055      & 583       & 9.70              & 14.19            \\
R       & E        & 3        & 35 & 902          & 512          & 847       & 673       & 6.49              & 23.87            \\
R       & E        & 4        & 35 & 923          & 455          & 921       & 582       & 0.22              & 21.79            \\
C       & T        & 3        & 35 & 543            & 638             & 604       & 547       & 5.70              & 10.02            \\
C       & T        & 4        & 35 & 781          & 412          & 665       & 455       & 17.44             & 9.45             \\
C       & E        & 3        & 35 & 608          & 445          & 585       & 532       & 3.93              & 16.29            \\
C       & E        & 4        & 35 & 730          & 388          & 700       & 480       & 4.29              & 19.23            \\
RC      & T        & 3        & 35 & 1368         & 594          & 1218      & 715       & 12.01             & 16.92            \\
RC      & T        & 4        & 35 & 1404         & 568          & 1288      & 633       & 9.04              & 10.25            \\
RC      & E        & 3        & 35 & 1135         & 537          & 1099      & 702       & 3.28              & 23.53            \\
RC      & E        & 4        & 35 & 1180         & 483          & 1160      & 547       & 1.68              & 11.70            \\ \hline
Avg.    &          &          &    &              &              &           &           & 6.79              & 16.06            \\
Max.    &          &          &    &              &              &           &           & 17.44             & 28.72            \\
Min.    &          &          &    &              &              &           &           & 0.22              & 7.39 \\
\hline
\end{tabular}
\end{table}

\begin{figure}[t!]
  \begin{subfigure}{0.33\textwidth}
    \includegraphics[width=\linewidth]{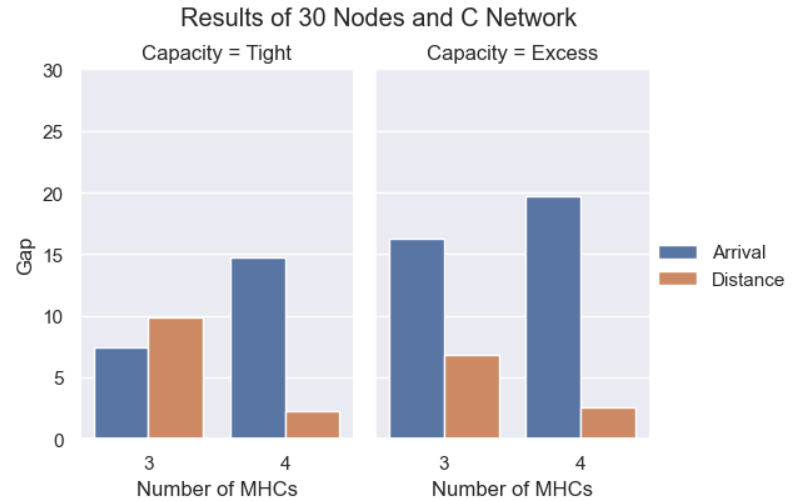}
    \caption{Bar plot of Network C}
    \label{fig:figure11}
  \end{subfigure}%
  \hfill
  \begin{subfigure}{0.33\textwidth}
    \includegraphics[width=\linewidth]{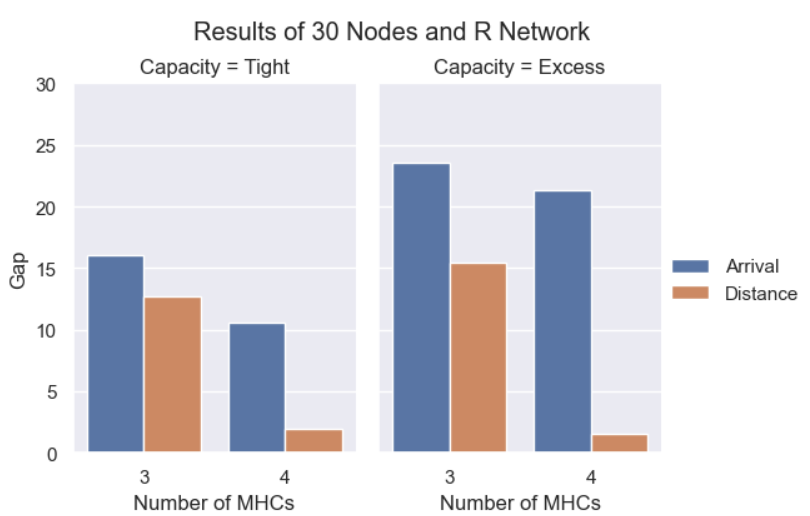}
    \caption{Bar plot of Network R}
    \label{fig:figure12}
  \end{subfigure}%
  \hfill
  \begin{subfigure}{0.33\textwidth}
    \includegraphics[width=\linewidth]{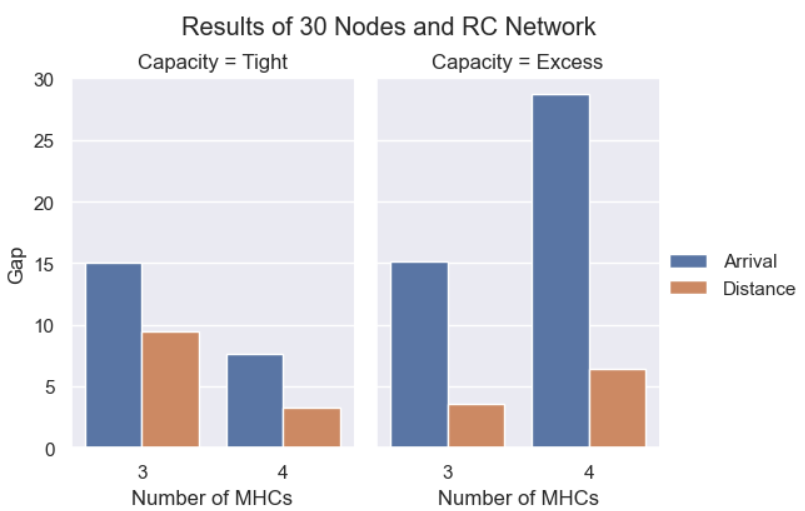}
    \caption{Bar plot of Network RC}
    \label{fig:figure13}
  \end{subfigure}
  \caption{Bar plot of 30 nodes.}
  \label{fig:bar30}
\end{figure}

\begin{figure}
  \begin{subfigure}{0.33\textwidth}
    \includegraphics[width=\linewidth]{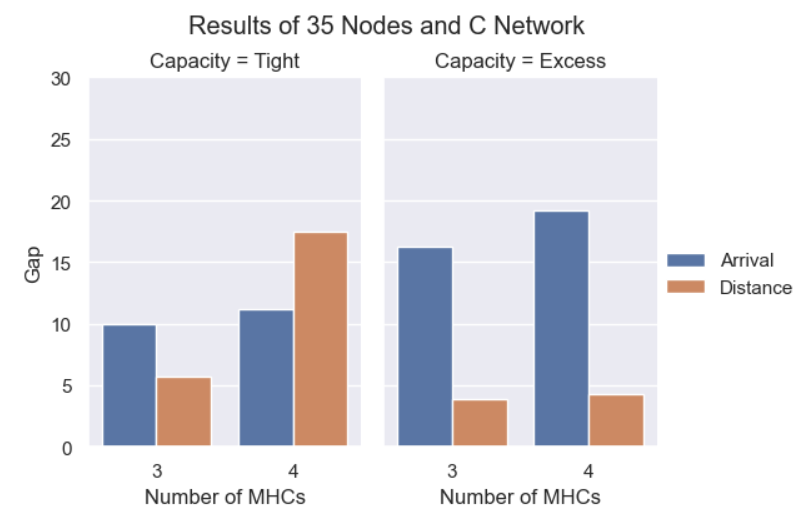}
    \caption{Bar plot of Network C}
    \label{fig:figure21}
  \end{subfigure}%
  \hfill
  \begin{subfigure}{0.33\textwidth}
    \includegraphics[width=\linewidth]{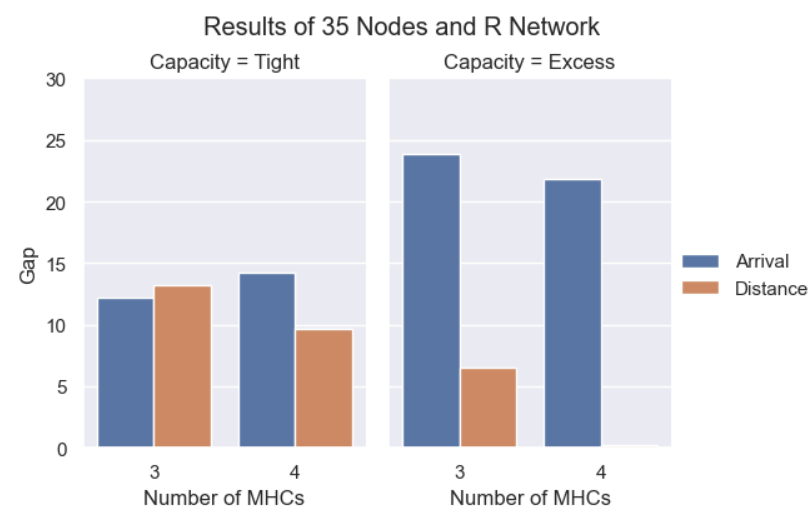}
    \caption{Bar plot of Network R}
    \label{fig:figure22}
  \end{subfigure}%
  \hfill
  \begin{subfigure}{0.33\textwidth}
    \includegraphics[width=\linewidth]{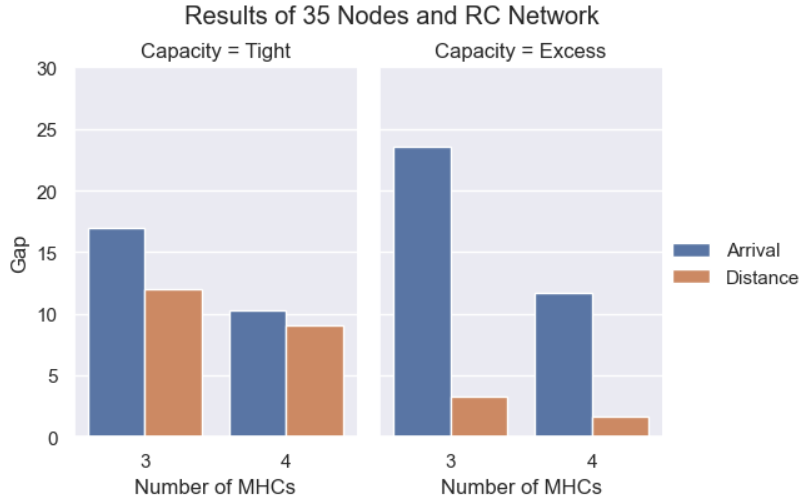}
    \caption{Bar plot of Network RC}
    \label{fig:figure23}
  \end{subfigure}
  \caption{Bar plot of 35 nodes.}
  \label{fig:bar35}
\end{figure}

We observe that the synchronization model is more expensive than the multi-trip model in terms of total traveled distance. This is reasonable given that the routing of the supply truck is a significant part of the total traveled distance. On the other hand, the synchronization model achieves a better outcome in terms of the arrival of the last MHC back to the depot. These results are important and valuable to organizations running MHC programs as it implies that by adding a truck to resupply the MHCs, the utilization of the MHCs improves significantly at a higher return considering that the total traveled distance metric. Figures \ref{fig:bar30} and \ref{fig:bar35} visualize these gains using grouped bar plots. We observe that it is better in all cases but two cases in which the increase is percentage of increase in traveled distance is higher than the reduction in the completion time of the last MHC arriving back to the depot. \\

\begin{figure}
  \begin{subfigure}{0.5\textwidth}
    \includegraphics[width=\linewidth]{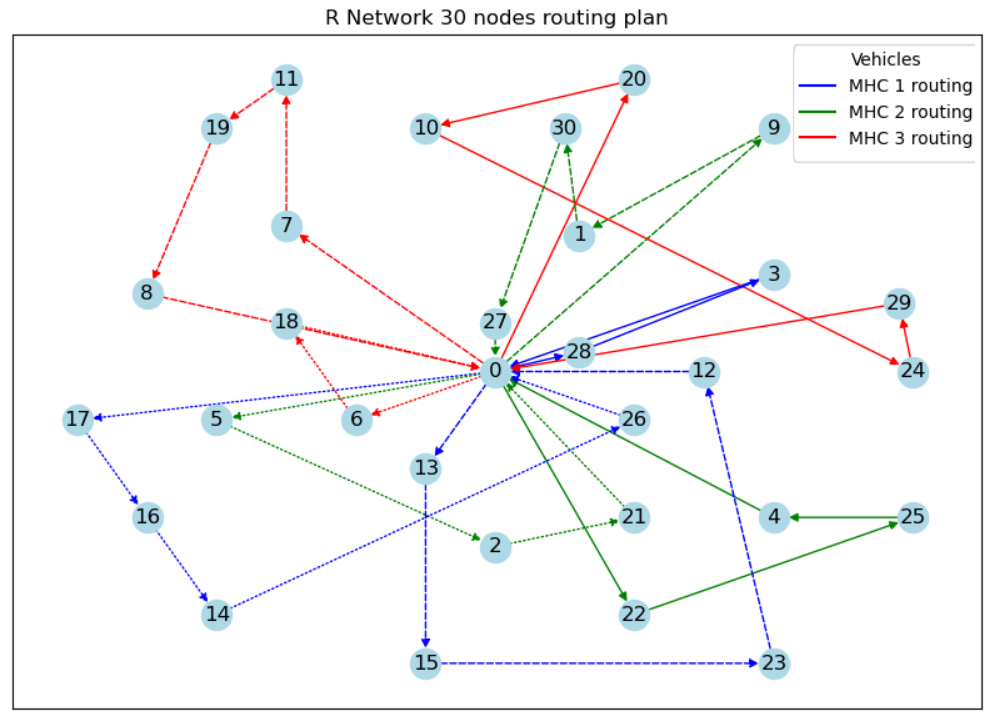}
    \caption{Routing plan under multi-trip model.}
    \label{fig:figure31}
  \end{subfigure}%
  \hfill
  \begin{subfigure}{0.5\textwidth}
    \includegraphics[width=\linewidth]{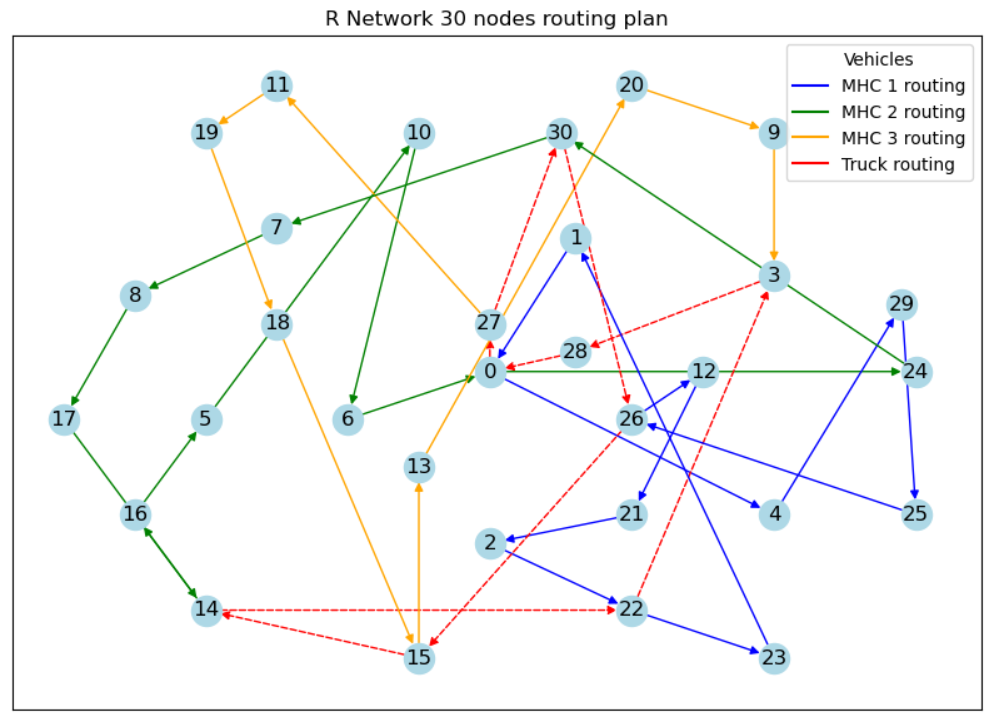}
    \caption{Routing plan under synchronization model.}
    \label{fig:figure32}
  \end{subfigure}%
  \caption{Routing plan of 30 nodes under tight MHC capacity.}
  \label{fig:routing30T}
\end{figure}

\begin{figure}
  \begin{subfigure}{0.5\textwidth}
    \includegraphics[width=\linewidth]{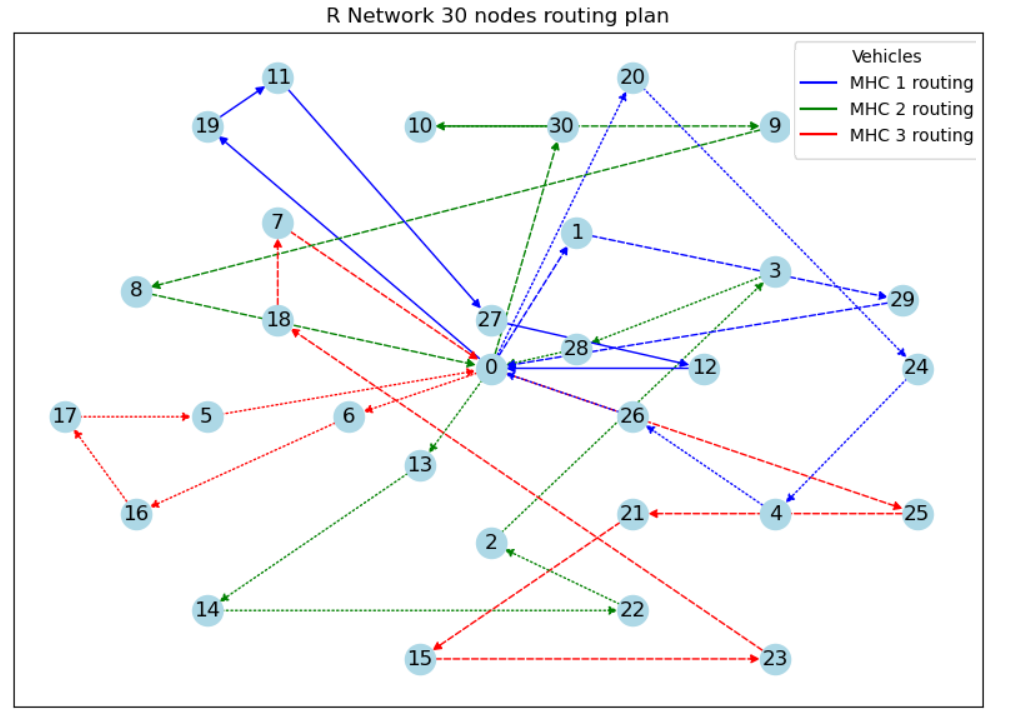}
    \caption{Routing plan under multi-trip model.}
    \label{fig:figure41}
  \end{subfigure}%
  \hfill
  \begin{subfigure}{0.5\textwidth}
    \includegraphics[width=\linewidth]{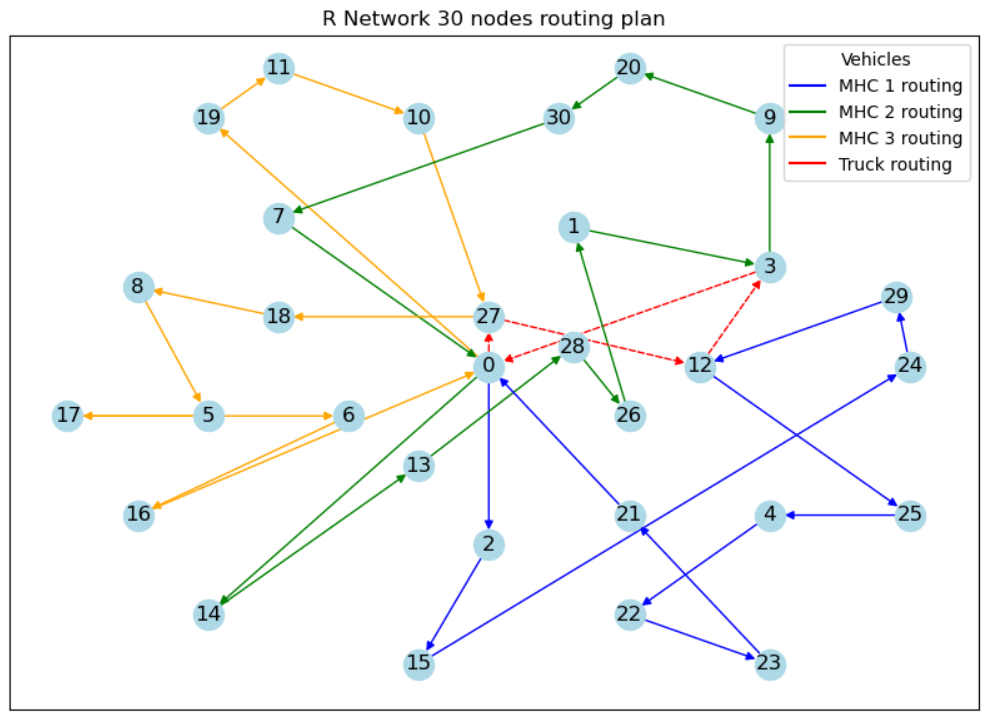}
    \caption{Routing plan under synchronization model.}
    \label{fig:figure42}
  \end{subfigure}%
  \caption{Routing plan of 30 nodes under excess MHC capacity.}
  \label{fig:routing30E}
\end{figure}

In Figures \ref{fig:routing30T} and \ref{fig:routing30E} we display the routing plan of the synchronization and the multi-trip models for a network with 30 nodes. For Figures \ref{fig:figure31} and \ref{fig:figure41}, colors indicate the routes of MHCs and the style of each line indicate the trip of the selected MHC.

\section{Conclusion}
\label{sec:conclude}
This paper investigated a new management model for the mobile health clinics where resupply process takes place en-route with a resupply truck. Instead of requiring the mobile health clinics to go back to the depot for a resupply, a truck is available to resupply an MHC en-route. A mixed integer programming model is developed, and Gurobi solver was used to solve small instances with 15-25 nodes. To solve larger instances of practical size, we proposed an adaptive large neighborhood search (ALNS) algorithm by combining classic and new destroy operators specially designed for the problem. Experimental tests on using Solomon's instances demonstrated that our proposed ALNS algorithm outperforms Gurobi solver on small instances.\\

To demonstrate the value of the synchronization approach, we compare the total distance traveled by the mobile health clinics and the resupply truck and arrival of the last mobile health clinic to the depot against a model where mobile health clinics are allowed to perform multiple trips to resupply from the depot. Our results reveal that despite the increase of 6.79\% in traveled distance under the synchronization approach, the reduction in the latest arrival is 16.06\% on average. Implying that the utilization of the fleet of mobile health clinics can be significantly improved at the cost of extra traveling time when combining the traveling time of all the mobile health clinics and the resupply truck. \\

In future work, our research will focus on incorporating uncertainity in service and traveling times. Adding uncertainity will be particularly challenging and valuable as the resource actions will be binary variables to coordinate the length of waiting time an MHC or truck may face. Another research extension that is worth studying is to extend the presented model to include planning horizon of multiple periods rather than a single period planning horizon model that we presented in this work.

\bibliographystyle{pomsref} 

 \let\oldbibliography\thebibliography
 \renewcommand{\thebibliography}[1]{%
 	\oldbibliography{#1}%
 	\baselineskip14pt 
 	\setlength{\itemsep}{10pt}
 }
\bibliography{main} 





\end{document}